\documentclass[12pt]{article}
\usepackage[paperwidth=230mm, paperheight=310mm]{geometry}
\usepackage{graphicx}
\usepackage{amssymb}
\usepackage{amsmath}
\usepackage{hyperref}
\usepackage{enumerate}
\usepackage{amsthm}
\usepackage{amsmath}
\usepackage{float}
\usepackage{mathtools}
 \usepackage{lineno}

\makeatletter
\renewenvironment{proof}[1][\proofname]{%
   \par\pushQED{\qed}\normalfont%
   \topsep6\p@\@plus6\p@\relax
   \trivlist\item[\hskip\labelsep\bfseries#1\@addpunct{.}]%
   \ignorespaces
}{%
   \popQED\endtrivlist\@endpefalse
}
\makeatother

\numberwithin{equation}{section}

\newtheorem{theorem}{Theorem}

 \numberwithin{theorem}{section}
 \newtheorem{corollary}[theorem]{Corollary}
\newtheorem{lemma}[theorem]{Lemma}

\newtheorem{remark}[theorem]{Remark}

\newtheorem{thmx}{Theorem}

\makeatletter
\def\keywords{\xdef\@thefnmark{}\@footnotetext}
\makeatother

\renewcommand{\P}{\mathbb{P}}
\newcommand{\E}{\mathbb{E}}
\newcommand{\R}{\mathbb{R}}

\newcommand{\cA}{\mathcal A}

\newcommand{\cD}{\mathcal D}

\newcommand{\cF}{\mathcal F}

\newcommand{\eps}{\varepsilon}

 \newcommand{\nn}{\nonumber}
 \newcommand{\no}{\noindent}

\newcommand{\Extra}[1]{{\color{blue}#1}}
\renewcommand{\Extra}[1]{}
\newcommand{\la}{\langle}
\newcommand{\ra}{\rangle}

\begin{document}

\author{
Jieliang Hong$^*$ \quad Leonid Mytnik$^\dagger$
}
\title{Exceptional times for the instantaneous propagation of superprocess}

\date{{\small  {\it  
$^{*}$Department of Mathematics, Southern University of Science and Technology, Shenzhen, China\\
$^{\dagger}$Faculty of Data and Decision Sciences, Technion, Haifa, Israel\\
 $^*$E-mail:  {\tt hongjl@sustech.edu.cn} \\
  $^\dagger$E-mail:  {\tt leonidm@technion.ac.il} 
  }
  }
  }
 
 \maketitle
\begin{abstract}
For a Dawson-Watanabe superprocess $X$ on $\R^d$, it is shown in Perkins \cite{Per90} that if the underlying spatial motion belongs to a certain class of L\'evy processes that admit jumps, then with probability one the closed support of $X_t$ is the whole space for almost all $t>0$ before extinction, the so-called ``instantaneous propagation'' property. In this paper for superprocesses on $\R^1$ whose spatial motion is the symmetric stable process of index $\alpha \in (0,2/3)$, we prove that there exist exceptional times at which the support is compact and nonempty. Moreover, we show that the set of exceptional times is dense with full Hausdorff dimension. Besides, we prove that near extinction, the support of the superprocess is concentrated arbitrarily close to the distinction point, thus upgrading the corresponding results in Tribe \cite{Tri92} from $\alpha \in (0,1/2)$ to $\alpha \in (0,2/3)$, and we further show that the set of such exceptional times also admits a full Hausdorff dimension.
\end{abstract}

\section{Introduction}

The Dawson-Watanabe superprocess is a measure-valued process that arises as a scaling limit of critical branching random walks. Let $M_F=M_F(\R^d)$ be the space of finite measures on $\R^d$ equipped with the topology of weak convergence of measures. Let $Y=\{Y_t, t\geq 0\}$ be a c\`adl\`ag Feller process on $\R^d$. Let $(\cA,\cD(\cA))$ be the generator of $Y$ where $\cD(\cA)$ is a suitable subset of bounded continuous functions.  A Dawson-Watanabe superprocess $X=(X_t,t\geq 0)$ starting from $X_0\in M_F$ is a continuous $M_F$-valued strong Markov process defined on some complete filtered probability space  $(\Omega, \cF, \cF_t, \P)$ such that $X$ satisfies the following martingale problem: 
\begin{align}\label{9e10.25}
(MP)_{X_0}^{\cA}: \quad &\text{ For any }  \phi \in \cD(\cA), \ M_t(\phi)=X_t(\phi)-X_0(\phi)-\int_0^t X_s(\cA\phi) ds\nn\\
&\text{ is  a continuous $(\cF_t)$-martingale with } \langle M(\phi)\rangle_t=\int_0^t X_s(\phi^2)ds.
\end{align}
Here and in what follows, for any $\mu \in M_F$ and integrable function $\phi$, we write 
\begin{align*}
\mu(\phi)=\int \phi(x) \mu(dx),
\end{align*}
 and $\mu(1)=\int 1   \mu(dx)$ is the total mass of $\mu$. By Theorem II.5.1 of \cite{Per02}, the above martingale problem $(MP)_{X_0}^{\cA}$ uniquely characterizes the law of $X$ (starting from $X_0$) on $C([0,\infty), M_F(\R^d))$, the space of continuous $M_F(\R^d)$-valued paths furnished with the compact-open topology. We denote such a law of the superprocess $X$ by $\P_{X_0}$.

Since the late eighties, the sample path properties of superprocesses have been studied extensively in, for example,  \cite{DIP89}, \cite{DP91}, \cite{Per88}, \cite{Per89}, \cite{Per90}. More recently, jointly with Edwin Perkins, we calculate in \cite{HMP18} the exact Hausdorff dimension of the topological boundary of the total support of the superprocess. 

Denote by $S(\mu)=\text{Supp}(\mu)$ the closed support of a measure $\mu$. One of the most impressive results for the sample path properties, called the ``instantaneous propagation'', concerns the support of $X$ when the underlying spatial motion is a L\'evy process that admits jumps.  

\begin{thmx}[Perkins \cite{Per90}]\label{0t1}
Let $X$ be a superprocess whose spatial motion is a L\'evy process with L\'evy measure $\nu$ satisfiying $\cup_{k=1}^\infty S(\nu^{(k)})=\R^d$, where $\nu^{(k)}$ is the $k$-fold convolution with itself. Then for any $X_0\in M_F(\R^d)$,  
\begin{align} 
\P_{X_0}(S(X_t)=\R^d|X_t\neq 0)=1, \ \forall t>0.
\end{align}
\end{thmx}
By applying Fubini's theorem, the above readily implies that with $\P_{X_0}$-probability one, $S(X_t)= \R^d$ for almost all $0<t<\zeta$, where $\zeta=\zeta_X=\inf\{t\geq 0: X_t=0\}$ is the extinction time of the superprocess. Similar results have been extended in Evans-Perkins \cite{EP91} to more general spatial motions under some conditions (see Corollary 5.3 of the same reference). Superprocesses with different branching structures exhibit the same propagation phenomena (see Li-Zhou \cite{LZ08}), meaning that the instantaneous propagation property is independent of the branching mechanism. A more recent paper Hughes-Zhou \cite{HZ22} proves that the $\Lambda$-Fleming-Viot process also shows similar instantaneous propagation properties. 

Such instantaneous propagation, however, does not describe the sample path of the superprocess precisely. Define the set of exceptional times to be
\begin{align} 
H=\{t>0: S(X_t) \text{ is compact and nonempty}\}.
\end{align}
It has been conjectured in Perkins \cite{Per90} that $H$ is nonempty a.s. for the symmetric stable superprocess with index $0<\alpha<2$, i.e. when the underlying spatial motion of the superprocess is a symmetric $\alpha$-stable process. Later Tribe studied in \cite{Tri92} the behavior of the superprocess near extinction and discovered that a sequence of exceptional times exists near extinction where the support is concentrated arbitrarily close to a single point, provided the generator $\cA$ of the spatial motion is bounded. In particular, Tribe \cite{Tri92} proves such results for the one-dimensional symmetric stable superprocess with index $0<\alpha<1/2$. 

For any $x\in \R^d$ and $r>0$, define 
\begin{align*}
B(x,r)=\{y\in \R^d: |y-x|<r\},   \text{ and set } B_r=B(0,r).
\end{align*}

\begin{thmx}[Tribe \cite{Tri92}] \label{0t2}
Let $X_0\in M_F(\R^d)$.

\noindent (i) For any superprocess $X$ satisfying $(MP)_{X_0}^{\cA}$, there exists a random variable $F$ on $\R^d$ such that with $\P_{X_0}$-probability one, 
\begin{align}\label{e4.71}
\frac{X_t}{X_t(1)}\to \delta_F \text{ in $M_F(\R^d)$ as $t\to \zeta$.} 
\end{align}
(ii) If the generator $\cA$ of the spatial motion is bounded or the spatial motion is a symmetric $\alpha$-stable process on $\R$ with $\alpha \in (0,1/2)$, then for all $r>0$, with $\P_{X_0}$-probability one, there exists $t_n\uparrow \zeta$ such that
\begin{align}\label{0e5.79}
S(X_{t_n}) \subseteq B(F,r),
\end{align}
where $F$ is as in \eqref{e4.71}.
\end{thmx}

 Tribe's theorem tells us that although the superprocess has the full support of the space almost all the time, the support collapses at some exceptional times near extinction, thus giving $H$ is nonempty.

Although the existence of exceptional times is settled for the case of bounded some generator and $\alpha \in (0,1/2)$, whether or not such support collapsing results in Theorem \ref{0t2} (ii) continue to hold for $\alpha \in [1/2, 2)$ has not be resolved since then. Moreover, the picture is still obscure for the sample path of the superprocess: For instance, we do not know if exceptional times exist elsewhere other than near extinction. It would be desirable to recognize whether the support-collapsing property of the superprocess is intrinsic or purely owing to that the population is fading out.

In this paper, we solve the above problems and demonstrate that exceptional times exist everywhere with full Hausdorff dimension. This provides a clearer picture of the sample path: Almost all the time, the support of the superprocess is the entire space; in the meantime, the support will collapse to a compact set instantly after any time.

 Denote by $\text{dim}(K)$ the Hausdorff dimension of any compact set $K$. 

\begin{theorem}\label{t1}
Let $X$ be a one-dimensional symmetric stable superprocess with index \break $\alpha \in (0,2/3)$ starting from $X_0\in M_F(\R)$. Then $\P_{X_0}$-a.s. that  $H$ is dense in $[0,\zeta]$ with $\text{dim}(H)=1$. Moreover, for any $t\geq 0$, 
\begin{align}\label{e0.0}
&\P_{X_0}\Big(\bigcap_{N=1}^\infty \Big\{\text{dim}\big(H \cap (t,t+N^{-1})\big)=1\Big\}\Big|X_t\neq 0\Big)=1.
\end{align}
 \end{theorem}

Next, somehow surprisingly, by using a similar proof to that of Theorem \ref{t1}, we extend Tribe's results for symmetric stable superprocesses in Theorem \ref{0t2} (ii) from $\alpha \in (0,1/2)$ to $\alpha \in (0,2/3)$. Moreover, instead of the simple existence of such exceptional times, we largely improve Tribe's result by showing that the set of times when the support collapses to the neighborhood of the extinction point has a full Hausdorff dimension.

For any $r>0$ and $x\in \R$, define
\begin{align} \label{0e5.99}
H^r_x=\Big\{t>0: S(X_t) \subseteq \overline{B(x,r)} \text{ and $S(X_t)$ is nonempty}\Big\}.
\end{align}

Note that in the above definition, we use the closed ball instead of the open ball in \eqref{0e5.79}. But of course, they are equivalent as $r>0$ is arbitrary.

 \begin{theorem}\label{t2}
Let $X$ be a one-dimensional symmetric stable superprocess with index \break $\alpha \in (0,2/3)$ starting from $X_0\in M_F(\R)$. For any $r>0$, with $\P_{X_0}$-probability one, for any $N\geq 1$, 
\[
 H^r_F \cap (\zeta-N^{-1},\zeta) \neq \emptyset, \quad \text{ and moreover,} \quad \text{dim}\Big(H^r_F \cap (\zeta-N^{-1},\zeta)\Big)=1,
\]
where $F$ is as in \eqref{e4.71}.
 
\end{theorem}

 \begin{remark}
 We are optimistic that the above two theorems shall hold in any dimensions, but we only prove $d=1$ in the current paper as the integral calculus is already quite involved for $d=1$; we leave this task for $d\geq 2$ to an energetic reader. The condition for  $\alpha \in (0,2/3)$ comes from some technical calculations for the moments of the immigration density (see, e.g., Lemma \ref{l1.33} below). We give some heuristic explanations of why requiring $\alpha \in (0,2/3)$ in Remark \ref{r4.2}.
 \end{remark}

The idea for the proof originated from the SPDE techniques used in \cite{MMQ11} (see also \cite{BMS21}) to study the front propagation speed for solutions to reaction-diffusion equations. Denote by $\Delta_\alpha=-(-\Delta)^{\alpha/2}$ the generator of the symmetric $\alpha$-stable process on $\R$.  
Consider the SPDE
\begin{align} \label{e3.1}
\frac{\partial X(t,x)}{\partial t} =\Delta_\alpha X(t,x)+\sqrt{X(t,x)} \dot{W}(t,x), \ \quad X\geq 0,
\end{align}
where $W=W(t,x)$ is a space-time white noise. 
 Set $R> 0$. Let $\Delta_\alpha^R$ be the generator of the symmetric $\alpha$-stable process {\it killed} when exiting $B_R$. Decompose the solution $X(t,x)$ from \eqref{e3.1} as $V^R_t+W^R_t$ where $V^R_t(x)=V^R(t,x)$ is a solution to the SPDE
\begin{align} \label{e3.2}
\begin{dcases}
\frac{\partial V^R_t(x)}{\partial t} =\Delta_\alpha^R V^R_t(x)+\sqrt{V^R_t(x)} \dot{W}(t,x), \quad &x\in B_R,\\
V^R_t(x)=0, \quad &x\in B_R^c.
\end{dcases}
\end{align}
By slightly abusing the notation, we use $S(f)$ to denote the support of a function $f$. Hence $S(V^R_t)\subseteq \overline{B_R}$. It suffices to show that there exist times at which $W^R_t(1)=0$ as it implies that $S(X_t)=S(V^R_t)\subseteq \overline{B_R}$ is compact. 

Alternatively, one may define $V^R=(V^R_t, t\geq 0)$ to be a superprocess whose underlying spatial motion is the $\alpha$-stable process killed at the exit of the ball $B_R$ so that the support of $V^R_t$ is contained in $\overline{B_R}$ for any $t>0$. Next, let $W^R=(W^R_t, t\geq 0)$ be the symmetric stable superprocess with extra immigration arising from those killed ``particles'' in the process $V^R$ so that $\{V^R_t+W^R_t, t\geq 0\}$ is equal in law to $\{X_t, t\geq 0\}$. For any $\delta>0$, by letting $R>0$ large or when the superprocess is near extinction, we will show that $\{W^R_t(1), t\geq 0\}$ is stochastically bounded above by the square of a $\delta$-dimensional Bessel process starting from $0$. Then the zeros of the Bessel processes will give the exceptional times at which the support of $X_t$ is in $\overline{B_R}$. The rigorous proof will be carried out in Section \ref{s2}. \\

\no {\bf Organization of the paper}.   In Section \ref{s2}, we give the proof of the main theorems assuming some moment results on the immigration term arising from the killed particles. In Section \ref{s3}, we use the moment formulas from Konno and Shiga \cite{KS88} to prove those moment estimates. The Appendix contains an elementary proof regarding the Bessel process and some technical calculus computations.\\

 \no ${\bf Notation\ and\ convention\ on\ constants.}$ Constants whose value is unimportant and may change from line to line are denoted $C, c, c_d, c_1,c_2,\dots$. In what follows, for positive integers $k\geq 1$, $C_b^k(\R)$ will denote the space of bounded continuous functions whose derivatives of order less than $k+1$ are also bounded continuous.  \\

\section*{Acknowledgements}
The authors' work was partly supported by ISF grant No. 1985/22. We thank Edwin Perkins for very helpful conversations and inspiring suggestions. We thank Zenghu Li for pointing out an error in the early version of the work.  We thank Renming Song for telling us the reference for the Dirichlet fractional Laplacian.

\section{Proof of the main results}\label{s2}

\subsection{Decomposition of superprocess} \label{s2.1}
Let $M_c(\R)$ be the space of compactly supported finite measures on $\R$ equipped with the topology of weak convergence of measures.
 In this subsection, we will decompose the one-dimensional symmetric stable superprocess $X=(X_t, t\geq 0)$,  with index $0<\alpha<2/3$, starting from some $X_0 \in M_c(\R)$. Recall that $X$ satisfies the martingale problem $(MP)_{X_0}^{\Delta_\alpha}$ as in \eqref{9e10.25}, that is, for any $f\in C_b^2(\R)$,
  \begin{align*} 
 X_t(f)= X_0(f)+M_t(f)+\int_0^t   X_s(\Delta_\alpha f)  ds,
 \end{align*}
where $\Delta_\alpha=-(-\Delta)^{\alpha/2}$ is the generator of the symmetric $\alpha$-stable  process on $\R$, also called the fractional Laplacian. The integral representation of $\Delta_\alpha$ is given by (see, e.g., Definition 2.5 of \cite{Kwa17})
   \begin{align}\label{4e7.80}
 \Delta_\alpha f(x)=\lim_{\eps\downarrow 0} \int_{\R-B(x,\eps)} \Big(f(y)-f(x)\Big) \frac{c_\alpha}{|y-x|^{1+\alpha}} dy.
 \end{align} 
 Here $c_{\alpha}>0$ is some constant depending on $\alpha$ (see \cite{Kwa17} for the explicit expression).

Assume that $X_0\in M_c(\R)$ and $R>0$ satisfy
 \begin{align}\label{e1.04}
 \text{Supp}(X_0)\subset \overline{B(0,R/2)}. 
 \end{align}
Now consider the symmetric $\alpha$-stable process killed when exiting ${B_R}$, the generator of which is called the Dirichlet fractional Laplacian $\Delta_\alpha^R:=(\Delta_\alpha)|_{B_R}$. It is well known that for $f\in C_b^2(\R)$,
   \begin{align}\label{4e7.81}
\Delta_\alpha^R f(x)=(\Delta_\alpha)|_{B_R}f(x)=\lim_{\eps\downarrow 0}  \int_{B_R-B(x,\eps)}& (f(y)-f(x))\frac{c_\alpha}{|y-x|^{1+\alpha}}dy\nn\\
&-f(x) \int_{B_R^c}\frac{c_\alpha}{|y-x|^{1+\alpha}}dy, \quad x\in B_R,
 \end{align}
 where $B_R^c$ is the complement of $B_R$.
We refer the reader to Theorem 4.4.3 of Fukushima-Oshima-Takeda \cite{FOT94} for the Dirichlet form of the killed $\alpha$-stable process, from which one may obtain \eqref{4e7.81} immediately. Combine \eqref{4e7.80} and \eqref{4e7.81} to get
   \begin{align}\label{5e11.333}
 (\Delta_\alpha- \Delta_\alpha^R) f(x):= \Delta_\alpha f(x)- \Delta_\alpha^R f(x)=& \int_{ B_R^c} f(y)  \frac{c_\alpha}{|y-x|^{1+\alpha}} dy,\quad x\in B_R.
 \end{align}
 
Consider a superprocess $V^R=(V_t^R, t\geq 0)$ starting from $V_0^R=X_0$, whose spatial motion is the killed $\alpha$-stable process generated by $\Delta_\alpha^R$, that is, $V^R$ a solution to the martingale problem $(MP)_{X_0}^{\Delta_\alpha^R}$ as in  \eqref{9e10.25}   so that for any $f\in C_b^2(\R)$,
 \begin{align}\label{5e3}
 V_t^R(f)=X_0(f)+M^{V^R}_t(f)+\int_0^t   V_s^R(\Delta_\alpha^R f)  ds,
 \end{align}
  where $M^{V^R}(f)$ is a continuous martingale whose quadratic variation is given by \break $t\mapsto \int_0^t  V^R_s(f^2) ds$.
  
  Next, for any $f\in C_b^2(\R)$, define
  \begin{align} \label{be2.1}
{A}_t^R(f)= \int_0^t  V_s^R(\Delta_\alpha f-\Delta_\alpha^R f ) ds =\int_0^t   ds \int V_s^R(dx) \int_{ B_R^c}  f(y) \frac{c_\alpha}{|y-x|^{1+\alpha}} dy.
 \end{align}
 Set 
 \begin{align} \label{4e2.4}
{A}_t^R:={A}_t^R(1) \text{ and define }    \dot{A}_t^R:=   \int V_t^R(dx) \int_{ B_R^c}   \frac{c_\alpha}{|y-x|^{1+\alpha}} dy.
 \end{align}
 
      \begin{lemma}\label{0l4.1} 
 For any $\alpha \in (0, 2/3)$, there exist some constant $C>0$ depending only on $\alpha$ such that for any $0<\eps_0< 1\wedge (\frac{1}{\alpha}-\frac{3}{2})$, if $X_0\in M_c(\R)$ and $R>0$ satisfy \eqref{e1.04}, then
 \begin{align}\label{be3.42}
 \E\Big((\dot{A}_{t+s}^R -\dot{A}_t^R )^{4}\Big) \leq   Cs^{1+\eps_0} (X_0(1) \vee X_0(1)^4)   (R^{-\alpha} \vee R^{-8\alpha}), \quad \forall 0\leq t,s \leq 1.
\end{align}
\end{lemma}  
 The proof of the above lemma is deferred to Section \ref{s3}.
By applying Kolmogorov's continuity criterion,  \eqref{be3.42} shows that ($\dot{A}_t^R, t\geq 0)$ admits a continuous version. Since $$|{A}_t^R(f)-{A}_s^R(f)|\leq \|f\|_\infty \int_s^t \dot{A}_r^Rdr, \quad \forall t\geq s\geq 0,$$ we conclude that ${A}_t^R(f)$ is well-defined for all  $f\in C_b^2(\R)$. Moreover, $t\mapsto  {A}_t^R(f)$ is absolutely continuous, and $\dot{A}_t^R$ defined in \eqref{4e2.4} is indeed the derivative of $A_t^R$ in $t$.\\

 Now let $W^R=(W_t^R, t\geq 0)$ be a superprocess such that for any $f\in C_b^2(\R)$,
 \begin{align}\label{5e2}
 W_t^R(f)= M^{W^R}_t(f)+\int_0^t  W_s^R(\Delta_\alpha f)  ds + {A}_t^R(f),
 \end{align}
 where $M^{W^R}(f)$ is a continuous martingale with quadratic variation $t\mapsto \int_0^t W^R_s(f^2) ds$, and is orthogonal to $M^{V^R}(f)$. The existence and uniqueness in law of such a superprocess $W^R$ should be trivial if we consider it as a limit of the branching particle system with immigration term ${A}_t^R(f)$ arising from those killed particles in $V^R$.
 Combine \eqref{5e3}, \eqref{be2.1} and \eqref{5e2} to get
  \begin{align*} 
 V_t^R(f)+W_t^R(f)= &X_0(f)+\big(M^{V^R}_t(f)+M^{W^R}_t(f)\big)\nn\\
 &+\Big[\int_0^t   W_s^R(\Delta_\alpha  f)  ds+\int_0^t   V_s^R(\Delta_\alpha  f)  ds\Big], \ \forall f\in C_b^2(\R), t\geq 0.
  \end{align*}
  One can check that the process $(V_t^R+W_t^R, t\geq 0)$ satisfies the martingale problem $(MP)_{X_0}^{\Delta_\alpha}$ as in \eqref{9e10.25}, so $(V_t^R+W_t^R, t\geq 0)$  equals in law to that of $(X_t, t\geq 0)$ by the uniqueness of the law for  $(MP)_{X_0}^{\Delta_\alpha}$ (see, e.g., Theorem II.5.1 of \cite{Per02}). Hence we may define $X, V^R, W^R$ on a common probability space $(\Omega, \cF, \cF_t, \P)$ such that
    \begin{align} \label{5e2.87}
 X_t(\cdot)=V_t^R(\cdot)+W_t^R(\cdot), \quad \forall t\geq 0.
  \end{align}

 Take $f\equiv 1$ in  \eqref{5e2} and then use \eqref{be2.1} and \eqref{4e2.4} to see that
  \begin{align} \label{5e2.456}
 W_t^R(1)&=  M^{W^R}_t(1)+{A}_t^R=M^{W^R}_t(1)+\int_0^t \dot{A}_s^Rds, \quad \forall t\geq 0.
 \end{align}
 Recall that $M^{W^R}(1)$ is a continuous martingale with quadratic variation $t\mapsto \int_0^t W^R_s(1) ds$. By the martingale representation theorem, there is some standard one-dimensional Brownian motion $\{\beta_t, t\geq 0\}$ such that
  \begin{align}\label{5e2.4}
 W_t^R(1)=  \int_0^t \sqrt{W^R_s(1)} d\beta_s +\int_0^t      \dot{A}_s^R  ds, \quad \forall t\geq 0.
 \end{align}

The above equation allows us to compare $W_t^R(1)$ with the square Bessel process. We are now ready to prove Theorem \ref{t1} and Theorem \ref{t2}.
  
  \subsection{Proof of Theorem \ref{t1}}

 Fix $0<\alpha<2/3$. Let $X$ be a one-dimensional symmetric stable superprocess with index $\alpha$ starting from some $X_0\in M_F(\R)$.  By a simple scaling argument, we may assume that $X_0(1)=1$ without loss of generality. For any $K>0$, set  
\begin{align}\label{e1.33}
X_0^{1,K}(\cdot)=X_0(\cdot \cap B_K) \quad \text{ and } \quad X_0^{2,K}(\cdot)=X_0(\cdot \cap B_K^c).
\end{align}
Let $X^{1,K}$ and $X^{2,K}$ be two independent symmetric stable superprocess with index $\alpha$ starting respectively from $X_0^{1,K}$ and $X_0^{2,K}$. Define $X,$ $X^{1,K}$ and $X^{2,K}$ on a common probability space $(\Omega, \cF, \cF_t, \P)$ such that
\[
X_t(\cdot)=X^{1,K}_t(\cdot)+X^{2,K}_t(\cdot), \quad \forall t\geq 0.
\]

 Fix any $\eps \in (0,1/4)$. Let $\zeta_{X^{i,K}}$ be the extinction time of $X^{i,K}$ for $i=1,2$. By (II.5.12) of \cite{Per02}, we get
\begin{align}\label{ea1.33}
\P(\zeta_{X^{i,K}}>t)=\P(X^{i,K}_t(1)>0)=1-\exp(- {2X^{i,K}_0(1)}/{t}), \quad \forall t> 0.
\end{align}
 Since $X_0(1)=1$, we may pick $K=K(\eps)>0$ to be large enough such that
 \[
 X^{2,K}_0(1)=X_0(B_K^c)<\eps^3/2  \quad  \text{ and } \quad  X^{1,K}_0(1)= 1- X^{2,K}_0(1)>1-\eps^3/2 >1/2. 
 \]
By using \eqref{ea1.33} and the above, we get
\[
\P(\zeta_{X^{1,K}}>\eps)=1-\exp(- {2X^{1,K}_0(1)}/{\eps})\geq 1-e^{-\eps^{-1}}\geq 1-\eps,
\]
and
 \[
\P(\zeta_{X^{2,K}}>\eps^2)=1-\exp(- {2X^{2,K}_0(1)}/{\eps^2})\leq 2 X^{2,K}_0(1)/\eps^2\leq   \eps.
\]
 Combine the above two inequalities to obtain
   \begin{align}\label{5e23.2}
\P(\zeta_{X^{2,K}}\leq \eps^2<\eps<\zeta_{X^{1,K}})\geq 1-2\eps. 
   \end{align}

Next,  since $X_0^{1, K}$ is supported on $\overline{B_K}$, we may let $R>2K+1$ so that \eqref{e1.04} is satisfied with $X_0=X_0^{1, K}$. Following the derivation of \eqref{5e2.87}, one may further define $V^R$ as in \eqref{5e3} with $X_0=X_0^{1, K}$ and  $W^R$ as in \eqref{5e2} such that
\[
X^{1,K}_t(\cdot)=V^{R}_t(\cdot)+W^{R}_t(\cdot), \quad \forall t\geq 0.
\]  
In particular, we recall from \eqref{5e2.4} that $W^R(1)$ satisfies that
   \begin{align}\label{be2.31}
 W_t^R(1)=  \int_0^t \sqrt{W^R_s(1)} d\beta_s +\int_0^t      \dot{A}_s^R  ds, \quad \forall t\geq 0.
 \end{align}
Here $\dot{A}_s^R$ is defined as in \eqref{4e2.4}. Moreover, by Lemma \ref{0l4.1}, we get
 \begin{align}\label{ce3.42}
 \E\Big((\dot{A}_{t+s}^R -\dot{A}_t^R )^{4}\Big) \leq   C   s^{1+\eps_0}  (X_0^{1, K}(1)\vee X_0^{1, K}(1)^4)(R^{-\alpha} \vee R^{-8\alpha}) , \quad \forall 0\leq t,s \leq 1,
\end{align}
where $\eps_0=\frac{1}{2}(1\wedge (\frac{1}{\alpha}-\frac{3}{2}))$ and $C>0$ depends only on $\alpha$.

Fix any $\delta \in (0,1/4)$. The following result is an easy consequence of \eqref{ce3.42}.

 \begin{corollary}\label{p1}
There exists some constant $R_0\geq 1$, depending only on $\alpha,  \eps,\delta$, such that for any $R>R_0$,
\[
\P(\dot{A}_s^R < \delta, \quad \forall 0\leq s\leq 1) \geq 1-\eps.
\]
\end{corollary}  
 
 \begin{proof} 
Use $X_0^{1, K}(1)\leq 1$ and $R\geq 1$ to see that \eqref{ce3.42}   becomes
 \begin{align}\label{ce4.37}
 \E\Big((\dot{A}_{t+s}^R -\dot{A}_t^R )^{4}\Big) \leq   C s^{1+\eps_0} R^{-\alpha}, \quad \forall 0\leq t,s \leq 1.
\end{align}
 By Kolmogorov's continuity criterion (see, e.g. Corollary 1.2 of Walsh \cite{Wal86}), there exist positive constants $C$ and $\gamma$, depending only on $\eps_0$, and a random variable $\xi$ such that with probability one, 
 \begin{align}\label{e4.37}
 |\dot{A}_{t}^R-\dot{A}_{s}^R|\leq \xi |t-s|^{\eps_0/4} \Big(\log \frac{\gamma}{|t-s|}\Big)^{1/2}, \forall 0\leq s,t\leq 1,
 \end{align} 
 and $\E(\xi^4)\leq CR^{-\alpha}$. Recall from \eqref{4e2.4} that 
    \begin{align*} 
  \dot{A}_0^R=   \int V_0^R(dx) \int_{B_R^c}   \frac{c_\alpha}{|y-x|^{1+\alpha}} dy\leq \int \frac{C}{(R-|x|)^\alpha} V_0^R(dx).
   \end{align*}
Use $S(V_0^R)=S(X_0^{1,K}) \subseteq \overline{B(0,R/2)}$  to see that
 \begin{align}\label{be6.1}
  \dot{A}_0^R\leq C (R/2)^{-\alpha} X_0^{1,K}(1)\leq CR^{-\alpha}.
   \end{align}
Pick $R>0$ large such that $ \dot{A}_{0}^R\leq \delta/2$. Together with \eqref{e4.37}, we get
 \begin{align}\label{e4.38}
 \sup_{0\leq t\leq 1} \dot{A}_{t}^R \leq \dot{A}_{0}^R+\xi  \sup_{t\in (0,1]}t^{\eps_0/4} \Big(\log \frac{\gamma}{t}\Big)^{1/2} \leq \delta/2+C_{\eps_0,\gamma}\xi.
 \end{align} 
  Through Markov's inequality, we have
 \begin{align}\label{e4.39}
  \P(C_{\eps_0,\gamma}\xi \geq \delta/2)\leq \Big(\frac{2C_{\eps_0,\gamma}}{\delta}\Big)^4 \E(\xi^4)\leq C R^{-\alpha}\leq \eps,
 \end{align} 
  if we let $R>0$ be large. Combine \eqref{e4.38} and  \eqref{e4.39} to obtain
   \begin{align*}
 \P\Big(\sup_{0\leq t\leq 1} \dot{A}_{t}^R\geq \delta\Big) \leq \P(C_{\eps_0,\gamma}\xi \geq \delta/2) \leq \eps,
 \end{align*} 
 as required.
 \end{proof}

 Fix $R>R_0$ with $R_0$ from Corollary \ref{p1} such that
   \begin{align}\label{4e23.1}
\P(\dot{A}_s^R < \delta, \quad \forall 0\leq s\leq 1) \geq 1-\eps.
 \end{align}
 Let $\{Z_t, t\geq 0\}$ be a $4\delta$-dimensional square Bessel process starting from $0$ satisfying (see, e.g. Chapter XI of \cite{RY94})
   \begin{align} \label{0e23.1}
 Z_t=  2\int_0^t  (Z_s)^{1/2} d\beta_s + 4 \delta t,  \quad \forall t\geq 0,
 \end{align}
 where $\{\beta_t, t\geq 0\}$ is the same Brownian motion that drives the equation \eqref{be2.31} for $W^R(1)$.
By letting $\widetilde{Z}_t=Z_t/4$ for all $t\geq 0$, we get
   \begin{align}\label{5e23}
 \widetilde{Z}_t=  \int_0^t   (\widetilde{Z}_s)^{1/2} d\beta_s +   \delta t, \quad \forall t\geq 0.
 \end{align}

  Define $\tau=\inf\{s\geq 0:  \dot{A}_s^R\geq \delta\}$. By following the classical comparison principle (see, e.g., Yamada \cite{Yam73}), one can easily show by \eqref{be2.31} and \eqref{5e23} that
 \[
 \text{ with probability one, } W_t^R(1)\leq  \widetilde{Z}_t \text{ holds for all }  0\leq t<\tau.
 \]
By \eqref{4e23.1},  we get
 \[
 \P(\tau>1)= \P(\dot{A}_s^R < \delta, \quad \forall 0\leq s\leq 1) \geq 1-\eps,
 \]
and so it follows that
   \begin{align}\label{4e3.23} 
 \P(W_t^R(1)\leq  \widetilde{Z}_t, \forall 0\leq t\leq 1)\geq  \P(\tau>1) \geq 1-\eps.
 \end{align}

 For any function $f: [0,\infty)\to \R$,   define
\begin{align*}
\text{Zeros}(f)=\{t\geq 0: f(t)=0\}
\end{align*}
to be the zero set of $f$.

\begin{lemma}\label{l2.2}
For any $\delta \in (0,1/4)$, if $Z$ is a $4\delta$-dimensional square Bessel process starting from $0$, then there exists some constant $C_\delta>0$ depending only on $\delta$ such that
\begin{align*}
&\P(\text{Zeros}(Z) \cap (a, b)\neq \emptyset)\geq 1-C_\delta a^{1-2\delta} (b-a)^{2\delta-1}, \quad \forall 0<a<b.
\end{align*}
\end{lemma}
The proof of Lemma \ref{l2.2} is quite elementary and is deferred to Appendix \ref{4a}.  Together with the classical results on the Hausdorff dimension of the zero sets of the $4\delta$-dimensional Bessel process (see, e.g., Theorem 4.2 of \cite{LX98}), we may further obtain 
\begin{align*}
&\P\Big(\dim(\text{Zeros}(Z) \cap (\eps^2, \eps))= 1-2\delta\Big)\geq \P(\text{Zeros}(Z) \cap (\eps^2, \eps)\neq \emptyset)\\
&\geq 1-C_\delta \eps^{2-4\delta} (\eps-\eps^2)^{2\delta-1} \geq  1-C_\delta \eps^{1/2},
\end{align*}
where in the second inequality we have used Lemma \ref{l2.2} with $a=\eps^2, b=\eps$, and the last inequality follows by $\eps<1/4$ and $\delta<1/4$.
Use the above with \eqref{4e3.23} to get (recall $\widetilde{Z}_t=Z_t/4$)
\begin{align*}
&\P\Big( \dim(\text{Zeros}(W^R(1)) \cap (\eps^2, \eps)) \geq 1-2\delta\Big) \\
&\geq  \P\Big(\Big\{W_t^R(1)\leq  \widetilde{Z}_t,  \forall 0\leq t\leq 1\Big\} \bigcap\Big\{\dim(\text{Zeros}(Z) \cap (\eps^2, \eps))= 1-2\delta\Big\} \Big)\\
& \geq 1-(\eps+C_\delta  \eps^{1/2}).
\end{align*}
 
Combining the above with \eqref{5e23.2}, we get that 
\begin{align*}
&  \P\Big(\{ \zeta_{X^{2,K}}\leq \eps^2 <\eps<\zeta_{X^{1,K}}\} \bigcap\Big\{\dim(\text{Zeros}(W^R(1)) \cap (\eps^2, \eps)) \geq 1-2\delta\Big\} \Big)\\
& \geq 1-(3\eps+C_\delta  \eps^{1/2}).
\end{align*}
Hence with probability greater than or equal to $1-(3\eps+C_\delta  \eps^{1/2})$, the following event holds:
\begin{align}\label{5e23.4}
\text{for any $t\in (\eps^2, \eps)$},& \text{ $X_t^{2,K}(1)=0$, $X_t^{1,K}(1)=V_t^R(1)+W_t^R(1)>0$, }\nn\\
& \text{and }\text{ $\dim(\text{Zeros}(W^R(1)) \cap (\eps^2, \eps))\geq 1-2\delta$}.
\end{align}
Recall $H$ is the set of exceptional times when $X$ is not zero and is compactly supported. One can easily check that on the event \eqref{5e23.4}, for each $t\in \text{Zeros}(W^R(1)) \cap (\eps^2, \eps)$,  
\[
X_t(\cdot)=V_t^R(\cdot) \text{ is compactly supported and $X_t(1)>0$}.
\]
Therefore on the event \eqref{5e23.4}, we get $\text{Zeros}(W^R(1)) \cap (\eps^2, \eps) \subseteq H$, thus giving 
\begin{align}\label{5e23.3}
\P(\dim(H\cap (\eps^2, \eps))\geq 1-2\delta)\geq 1-(3\eps+C_\delta  \eps^{1/2}).
\end{align}
 For each $N\geq 1$, pick $\eps<N^{-1}$ such that
\begin{align*} 
\P\Big( \dim(H \cap (0,N^{-1}))\geq  1-2\delta\Big)\geq \P(\dim(H\cap (\eps^2, \eps))\geq 1-2\delta) \geq 1-(3\eps+C_\delta    \eps^{1/2}).
\end{align*}
Let $\eps\downarrow 0$ to get $\P(\dim(H \cap (0,N^{-1}))\geq  1-2\delta)=1$, and then let $\delta\downarrow 0$ to get $$\P(\dim(H \cap (0,N^{-1}))= 1)=1$$ holds for all $N\geq 1$.
Hence 
\[
\P\Big(\bigcap_{N=1}^\infty \Big\{\text{dim}\big(H \cap (0,N^{-1})\big)=1\Big\}\Big)=1,
\]
and \eqref{e0.0} for $t=0$ follows. The case for \eqref{e0.0} with $t>0$ follows easily by the Markov property.
 
Notice that for each $t\geq 0$, on the event $$\bigcap_{N=1}^\infty \Big\{\dim\big(H \cap (t,t+N^{-1})\big)= 1\Big\},$$ we get $t$ is a limiting point of ${H}$, so by \eqref{e0.0},
\[
\P\Big(\text{$t$ is a limiting point of ${H}$}\Big|X_t\neq 0\Big)=1, \forall t\geq 0.
\]
By Fubini's theorem, the above implies that with $\P_{X_0}$-probability one,  for almost all $0\leq t<\zeta$, $t$ is a limiting point of ${H}$, thus giving $H$ is dense in $[0,\zeta]$.  The proof of Theorem \ref{t1} is now complete.

 \subsection{Proof of Theorem \ref{t2}}
 
Fix $0<\alpha<2/3$, $X_0\in M_F(\R)$ and $r>0$. Let $\{x_n, n\geq 1\}$ be a countable collection of distinct points on $\R$ such that 
\begin{align}\label{e7.44}
\R=\bigcup_{n=1}^\infty B(x_n, r/4).
\end{align}
 Fix such a collection $\{x_n\}$ throughout the rest of this subsection. For each $y\in \R$, define a bounded continuous function $\phi_y^r$ such that 
 \begin{align}\label{0e5.89}
\phi_y^r(x)=
\begin{cases}
1,& x\in B(y,r/2)^c;\\
0,&  x\in B(y,r/4);\\
\text{linear},&\text{ otherwise}.
\end{cases}
\end{align}

 For any $\eps \in (0,1/4)$, define a stopping time $\tau_{\eps}$ by
\begin{align} \label{0e5.98}
\tau_{\eps}:=\inf \Big\{t>0: 0<X_t(1) \leq \eps \text{ and  $\exists k\geq 1$ s.t. } \frac{X_t(\phi_{x_k}^r)}{X_t(1)}\leq \eps^3\Big\}.
\end{align}
Recall the extinction point $F$ from Theorem \ref{0t2} such that $\frac{X_t}{X_t(1)}$ converges weakly to $\delta_F$. 
Notice that by \eqref{e7.44}, with probability one there is some $k_0=k_0(\omega)\geq 1$ such that \break $F\in B(x_{k_0}, r/4)$ and so by \eqref{0e5.89}, we have 
 \begin{align} \label{0e5.88}
  \lim_{ t\uparrow \zeta}\frac{X_t(\phi_{x_{k_0}}^r)}{X_t(1)} = \delta_F(\phi_{x_{k_0}}^r)=0. 
\end{align}
Together with $\{X_t(1),t\geq 0\}$ is a continuous Feller's diffusion (see, e.g., Theorem II.1.2 of \cite{Per02}), we conclude from the above that with probability one $\tau_{\eps}<\zeta$, and there is some ${k_{\tau_{\eps}}}\geq 1$ such that
\begin{align}\label{0e5.12}
0<X_{\tau_{\eps}}(1) \leq  \eps \text{ and }   \frac{X_{\tau_{\eps}}(\phi_{x_{k_{\tau_{\eps}}}}^r)}{X_{\tau_{\eps}}(1)} \leq \eps^3.
\end{align}
Note that there might be more than one such ${k_{\tau_{\eps}}}\geq 1$ satisfying the above; if so, we simply pick the smallest $k\geq 1$ to be $k_{\tau_{\eps}}$.

Recall from \eqref{0e5.99} that for each $x\in \R$,
\begin{align}\label{0e5.92}
H_{x}^r=\Big\{t>0: \text{ $S(X_t)$ is nonempty and $S(X_t) \subseteq \overline{B({x},r)}$}\Big\}.
\end{align}
We first give the following result in light of \eqref{0e5.12}.
 
\begin{lemma}\label{al4.1}
 For any $r>0$ and $\eps', \delta \in (0,1/4)$, there exists some constant $\gamma_0\in (0,1/4)$ depending only on $r, \delta, \alpha$ and $\eps'$ such that for all $\eps<\gamma_0$, if $\mu \in M_F(\R)$ and $y\in \R$ satisfy
\begin{align}\label{0e5.13}
0<\mu(1) \leq \eps \text{ and }  \frac{\mu(\phi_y^r)}{\mu(1)}\leq \eps^3,
\end{align}
then 
\begin{align} \label{9e5.13}
 \P_{\mu} \Big(\text{dim}( H_y^r\cap (0, \infty)) \geq 1-2\delta \Big)\geq 1-\eps'.
\end{align}
\end{lemma}

Assuming the above lemma, we may finish the proof of Theorem \ref{t2}.

\begin{proof}[Proof of Theorem \ref{t2} assuming Lemma \ref{al4.1}]

Define
\begin{align}\label{0e5.95}
H^r=\Big\{t>0: \text{ $S(X_t)$ is nonempty and $\exists k\geq 1$ s.t. $S(X_t) \subseteq \overline{B({x_k},r)}$}\Big\}
\end{align}
so that
\begin{align} 
H^r=\bigcup_{k=1}^\infty  H_{x_k}^r.
\end{align}

Set $\delta \in (0,1/4)$ and  $\eps \in (0,1/4)$. Recall $\tau_{\eps}$ from \eqref{0e5.98} and $x_{k_{\tau_\eps}}$ from \eqref{0e5.12} to see that
\begin{align}\label{0e5.15}
&\P_{X_0}\Big(\text{dim}(H^r\cap (\tau_{\eps}, \infty)) \geq 1-2\delta \Big)\nn\\
&\geq \P_{X_0}\Big(\text{dim}(H_{x_{k_{\tau_\eps}}}^r\cap (\tau_{\eps}, \infty)) \geq 1-2\delta\Big)\nn\\
 &=\P_{X_0}\Big[\P_{X_0}\Big(\text{dim}(H_{x_{k_{\tau_\eps}}}^r\cap (\tau_{\eps}, \infty)) \geq 1-2\delta\Big|\cF_{\tau_\eps} \Big)\Big]\nn\\
 &= \P_{X_0}\Big[\P_{X_{\tau_{\eps}}}\Big(\text{dim}(H_{x_{k_{\tau_\eps}}}^r\cap (0, \infty)) \geq 1-2\delta \Big)\Big],
\end{align}
where the last equality follows by the strong Markov property of the superprocess.

 For any $\eps'>0$, let $\gamma_0 \in (0,1/4)$ be as in Lemma \ref{al4.1}. Then for any $\eps<\gamma_0$, we may use \eqref{0e5.12} and apply Lemma \ref{al4.1} with $\mu=X_{\tau_{\eps}}$ and $y=x_{k_{\tau_\eps}}$ to get
\begin{align*} 
\P_{X_{\tau_{\eps}}}\Big(\text{dim}(H_{x_{k_{\tau_\eps}}}^r\cap (0, \infty)) \geq 1-2\delta \Big) \geq 1-\eps'.
\end{align*}
 Hence \eqref{0e5.15} becomes
 \begin{align} \label{0e7.1}
&\P_{X_0}\Big(\text{dim}(H^r\cap (\tau_{\eps}, \infty)) \geq 1-2\delta\Big)\geq 1-\eps'.
\end{align}
Notice that with $\P_{X_0}$-probability one, $\tau_\eps\uparrow \zeta$ as $\eps\downarrow 0$. Hence for any $N\geq 1$, there exist some $\gamma_1\in (0,1/4)$ depending only on $X_0$, $\alpha, \eps'$ and $N$ such that for all $\eps<\gamma_1$,
 \begin{align} \label{0e7.2}
&\P_{X_0}\Big(\zeta-N^{-1}< \tau_{\eps}<\zeta\Big)\geq 1-\eps'.
\end{align}
For each $N\geq 1$, define the event
\[
A_N^{\delta}= \Big\{\text{dim}\Big(H^r\cap (\zeta-N^{-1}, \zeta)\Big) \geq 1-2\delta\Big\}.
\]
Combine \eqref{0e7.1} and \eqref{0e7.2} to conclude that for any $N\geq 1$, if $\eps<\gamma_0 \wedge \gamma_1$, then (recall $H^r\subseteq (0,\zeta)$)
 \begin{align*}  
 \P_{X_0}(A_N^{\delta})&\geq \P_{X_0}\Big(\Big\{\text{dim}\big(H^r\cap (\tau_{\eps}, \infty)\big) \geq 1-2\delta\Big\} \bigcap \Big\{\zeta-N^{-1}< \tau_{\eps}<\zeta \Big\}\Big) \\
 &\geq 1-2\eps'. 
\end{align*}
Since $\eps'>0$ is arbitrary and $A_N^{\delta}$ does not depend on $\eps'$, we may let $\eps' \downarrow 0$ to obtain $\P_{X_0}(A_N^{\delta})=1$. Define
\[
A_N= \Big\{\text{dim}\Big(H^r\cap (\zeta-N^{-1}, \zeta)\Big) =1\Big\} \text{ for each }N\geq 1.
\]
Let $\delta \downarrow 0$ to further get that $\P_{X_0}(A_N)=1$,
thus giving
 \begin{align} 
&\P_{X_0}\Big(\bigcap_{N=1}^\infty A_N\Big)=1.
\end{align}

Fix $\omega$ outside a null set so that $\cap_{N=1}^\infty A_N(\omega)$ holds. Recall from \eqref{e7.44} to see that there is some $k_0=k_0(\omega)\geq 1$ such that
 $F(\omega) \in B(x_{k_0(\omega)}, r/4)$ and (recall \eqref{0e5.88})
 \begin{align} \label{0e5.97}
 \frac{X_t(\phi_{x_{k_0}}^r)}{X_t(1)} \to 0 \text{ as $t\to \zeta$. }
\end{align}
It follows that there is some $N_0(\omega)\geq 1$ such that for all $N\geq N_0(\omega)$, we have  
 \begin{align*} 
 \frac{X_t(B(x_{k_0}, r/2)^c)}{X_t(1)} \leq  \frac{X_t(\phi_{x_{k_0}}^r)}{X_t(1)} \leq \frac{1}{4}, \quad \forall t\in (\zeta-N^{-1}, \zeta),
\end{align*}
where the first inequality is due to \eqref{0e5.89}, thus giving
 \begin{align} \label{0e5.77}
X_t\big({B(x_{k_0}, r/2)}\big)>0 , \quad \forall t\in (\zeta-N^{-1}, \zeta),
\end{align}

Fix any $N\geq N_0(\omega)$. For each $t\in H^r \cap (\zeta-N^{-1}, \zeta)$, the definition of $H^r$ from \eqref{0e5.95} implies that  there is some $k_1={k_1}(\omega)\geq 1$ such that $S(X_t) \subseteq \overline{B(x_{k_1}, r)}$. Together with \eqref{0e5.77}, it follows that 
 \begin{align}
 \overline{B(x_{k_0}, r/2)} \cap \overline{B(x_{k_1}, r)} \neq \emptyset,
 \end{align}
 and so $|x_{k_0}-x_{k_1}|\leq 3r/2$. Recalling that $F \in B(x_{k_0}, r/4)$, we obtain
 \begin{align*}  
S(X_t) \subseteq \overline{B(x_{k_1}, r)} \subseteq \overline{B(F,3r)}, \text{ and hence } t\in H_F^{3r}.
\end{align*}

 To make a summary, the above arguments show that with $\P_{X_0}$-probability one, for any $N\geq N_0(\omega)$, 
 \begin{align*}  
 H^r \cap (\zeta-N^{-1}, \zeta) \subseteq H_F^{3r} \cap (\zeta-N^{-1}, \zeta),
 \end{align*}
  thus giving $\text{dim}(H_F^{3r} \cap (\zeta-N^{-1}, \zeta) ) \geq \text{dim}(H^r \cap (\zeta-N^{-1}, \zeta))=1$. Therefore we conclude that for any $N\geq 1$,
\[
\P_{X_0}\Big(\text{dim}\Big(H^{3r}_F \cap (\zeta-N^{-1},\zeta)\Big)=1\Big)=1.
\]
It is immediate that 
\[
\P_{X_0}\Big( \bigcap_{N\geq 1} \Big\{\text{dim}\Big(H^{3r}_F \cap (\zeta-N^{-1},\zeta)\Big)=1\Big\}\Big)=1.
\]
Since $r>0$ is arbitrary, the proof of Theorem \ref{t2} is now complete.
\end{proof}

It remains to prove Lemma \ref{al4.1}.

\begin{proof} [Proof of Lemma \ref{al4.1}]
 Fix $r>0$. Let $X$ be a symmetric stable superprocess with index $\alpha\in (0,2/3)$ starting from some $\mu \in M_F(\R)$. 
By the translation invariance of superprocesses, we may assume that $\mu\in M_F(\R)$ satisfies \eqref{0e5.13} with $y=0$ and $\eps \in (0,1/4)$.
Set
\begin{align} \label{ce5.44}
\mu^{1,r}(\cdot)=\mu(\cdot \cap B(0,r/2)) \quad \text{ and } \quad \mu^{2,r}(\cdot)=\mu(\cdot \cap B(0,r/2)^c).
\end{align}
 Let $X^{1,r}$ and $X^{2,r}$ be two independent symmetric stable superprocess starting respectively from $\mu^{1,r}$ and $\mu^{2,r}$. Define $X,$ $X^{1,r}$ and $X^{2,r}$ on a common probability space $(\Omega, \cF, \cF_t, \P)$ such that
\begin{align*}
X_t(\cdot)=X^{1,r}_t(\cdot)+X^{2,r}_t(\cdot), \quad \forall t\geq 0.
\end{align*}

 Let $\zeta_{X^{i,r}}$ be the extinction time of $X^{i,r}$ for $i=1,2$. By \eqref{ea1.33}, we get for $i=1,2$,
\begin{align}\label{0a1.33}
\P(\zeta_{X^{i,r}}>t)=1-\exp(- {2\mu^{i,r}(1)}/{t}), \quad \forall t> 0.
\end{align}
Notice that by \eqref{ce5.44}, we get
\begin{align}\label{0a4.43}
\mu^{2,r}(1)=\mu(B(0,r/2)^c)\leq \mu(\phi_0^{r})\leq \mu(1) {\eps}^3,
\end{align}
where the first inequality uses \eqref{0e5.89} and the last inequality follows by the assumption \eqref{0e5.13} with $y=0$. 
Since $\eps^3<1/2$, we further obtain
\begin{align}\label{0a4.33}
\mu^{1,r}(1)=\mu(1)-\mu^{2,r}(1)\geq \frac{1}{2} \mu(1).
\end{align}
Set $i=1$ and $t=\eps \mu(1)$ in \eqref{0a1.33} to get that
\begin{align} 
\P(\zeta_{X^{1,r}}> \eps \mu(1))=1-\exp\Big(- \frac{2\mu^{1,r}(1)}{\eps\mu(1)}\Big) \geq 1-e^{-\eps^{-1}} \geq 1-\eps,
\end{align}
where the first inequality uses \eqref{0a4.33}. 
 For $\zeta_{X^{2,r}}$, we apply \eqref{0a1.33} with $i=2$ and $t=\eps^2 \mu(1)$ to obtain  
 \begin{align} 
\P(\zeta_{X^{2,r}}> \eps^2 \mu(1))=1-\exp\Big(- \frac{2\mu^{2,r}(1)}{\eps^2\mu(1)}\Big) \leq \frac{2\mu^{2,r}(1)}{\eps^2\mu(1)} \leq 2\eps,
\end{align}
where we have used \eqref{0a4.43} in the last inequality.  Combine the above two inequalities to conclude that
\begin{align}\label{0e2.77}
 \P(\zeta_{X^{1,r}}<\eps^2 \mu(1)<\eps \mu(1)<\zeta_{X^{1,r}}) \geq 1-3\eps.
\end{align}

Next,  since $\mu^{1, r}$ is supported on $\overline{B_{r/2}}$, we have \eqref{e1.04} is satisfied with $X_0=\mu^{1, r}$ and $R=r$. Following the derivation of \eqref{5e2.87}, one may further define $V^r$ as in \eqref{5e3} with $X_0=\mu^{1, r}$ and  $W^r$ as in \eqref{5e2} such that
\[
X^{1,r}_t(\cdot)=V^{r}_t(\cdot)+W^{r}_t(\cdot), \quad \forall t\geq 0.
\]  
In particular, we recall from \eqref{5e2.4} that $W^r(1)$ satisfies that
   \begin{align*} 
 W_t^r(1)=  \int_0^t \sqrt{W^r_s(1)} d\beta_s +\int_0^t      \dot{A}_s^r  ds, \quad \forall t\geq 0,
 \end{align*}
where $\dot{A}_s^r$ is defined by \eqref{4e2.4}. Moreover, Lemma \ref{0l4.1} gives that
 \begin{align}\label{de3.42}
 \E\Big((\dot{A}_{t+s}^r -\dot{A}_t^r )^{4}\Big) \leq   C  s^{1+\eps_0}  ( \mu^{1, r}(1)\vee  \mu^{1, r}(1)^4)(r^{-\alpha} \vee r^{-8\alpha}), \quad \forall 0\leq t,s \leq 1.
\end{align}
In the above, $\eps_0=\frac{1}{2}(1\wedge (\frac{1}{\alpha}-\frac{3}{2}))$ and $C>0$ depends only on $\alpha$.\\

Fix any $\eps', \delta \in (0,1/4)$. Recall that $\mu^{1, r}(1)\leq \mu(1) \leq \eps$. The following result is an easy consequence of \eqref{de3.42}.
 
 \begin{corollary} \label{0c3.1}
There exists some constant $\gamma_2 \in (0,1/4)$, depending only on $\alpha, r, {\eps'},\delta$, such that if $\eps\leq \gamma_2$, then 
\begin{align}\label{0e7.90}
\P(\dot{A}_s^r < \delta, \quad \forall 0\leq s\leq 1) \geq 1-{\eps'}.
\end{align}
\end{corollary}  
\begin{proof}
By using \eqref{de3.42} with $\mu^{1, r}(1)\leq  \eps<1$, we get
 \begin{align} 
 \E\Big((\dot{A}_{t+s}^r -\dot{A}_t^r )^{4}\Big) \leq   C(r) \eps s^{1+\eps_0}  , \quad \forall 0\leq t,s \leq 1,
\end{align}
The rest of the proof follows similarly to that of Corollary \ref{p1} using \eqref{ce4.37}. Note that in the proof here, $r>0$ is fixed and we set $\mu^{1, r}(1)\leq \eps<0$ to be small.
\end{proof}

Let $\eps \leq \gamma_2$ where $\gamma_2 \in (0,1/4)$ is as in Corollary \ref{0c3.1} so that \eqref{0e7.90} holds. Repeat the arguments from \eqref{4e23.1} to \eqref{4e3.23} to obtain
\begin{align}\label{9e4.44}
 \P(W_t^r(1)\leq  \widetilde{Z}_t, \forall  0\leq t\leq 1)  \geq  \P(\dot{A}_s^r < \delta,   \forall 0\leq s\leq 1)\geq 1-\eps',
 \end{align}
 where $\widetilde{Z}_t=Z_t/4$ and  $\{Z_t, t\geq 0\}$ is a $4\delta$-dimensional square Bessel process starting from $0$.  Apply Lemma \ref{l2.2} with $a=\eps^2 \mu(1)$ and $b=\eps \mu(1)$ to get
\begin{align}\label{9e3.44}
&\P(\dim(\text{Zeros}(Z) \cap (\eps^2 \mu(1),   \eps \mu(1)))= 1-2\delta) \nn\\
&\geq \P(\text{Zeros}(Z) \cap (\eps^2 \mu(1),   \eps \mu(1))\neq \emptyset)\geq 1-C_\delta \eps^{1-2\delta} \geq 1-C_\delta {\eps}^{1/2},
\end{align}
where the last  inequality follows by $\eps, \delta\in (0,1/4)$.
Combine \eqref{9e4.44} and \eqref{9e3.44} to get (recall $\widetilde{Z}_t=Z_t/4$)
\begin{align*}
&\P\Big( \dim(\text{Zeros}(W^r(1)) \cap (\eps^2 \mu(1),   \eps \mu(1))) \geq  1-2\delta\Big) \\
&\geq  \P\Big(\Big\{W_t^r(1)\leq  \widetilde{Z}_t,  \forall 0\leq t\leq 1\Big\}  \bigcap\Big\{\dim(\text{Zeros}(Z) \cap (\eps^2 \mu(1),   \eps \mu(1)))= 1-2\delta\Big\} \Big)\\
& \geq 1-(\eps'+C_\delta  \eps^{1/2}).
\end{align*}
Use the above with \eqref{0e2.77} to see that
\begin{align}\label{9e1.44}
& \P\Big(\Big\{\dim(\text{Zeros}(W^r(1)) \cap (\eps^2 \mu(1),   \eps \mu(1))) \geq  1-2\delta\Big\}\nn\\
&\quad \quad\quad \bigcap \Big\{\zeta_{X^{2,r}}<\eps^2 \mu(1)<\eps \mu(1)<\zeta_{X^{1,r}} \Big\}\Big) \nn\\
& \geq 1-(\eps'+C_\delta \eps^{1/2}+3\eps)\geq 1-5\eps',
\end{align}
where the last inequality follows by letting $\eps\leq \eps' \wedge (\eps'/C_\delta)^2$.\\

To make a summary, if $\eps \leq \gamma_2\wedge \eps' \wedge (\eps'/C_\delta)^2$, then with probability greater than or equal to $1-5\eps'$, the following event holds:
\begin{align}\label{0e23.4}
\text{for any $t\in (\eps^2 \mu(1),   \eps \mu(1))$},& \text{ $X_t^{2,r}(1)=0$, $X_t^{1,r}(1)=V_t^r(1)+W_t^r(1)>0$, }\nn\\
& \text{and }\text{ $\dim(\text{Zeros}(W^r(1)) \cap (\eps^2 \mu(1),   \eps \mu(1))\geq 1-2\delta$}.
\end{align}
One can easily check that on the event \eqref{0e23.4}, for each $t\in \text{Zeros}(W^r(1)) \cap (\eps^2 \mu(1),   \eps \mu(1))$,  
\[
X_t(\cdot)=X^{1,r}_t(\cdot)=V_t^r(\cdot) \text{ is supported on $\overline{B(0,r)}$ and $X_t(1)>0$}.
\]
Hence we recall $H_0^{r}$ from \eqref{0e5.92} to conclude that on the event \eqref{0e23.4},  
\[
\text{Zeros}(W^r(1)) \cap (\eps^2 \mu(1),   \eps \mu(1)) \subseteq H_0^{r}, 
\]
thus giving
\begin{align} 
\P(\dim( H_0^{r})\geq 1-2\delta)\geq 1-5\eps'.
\end{align}
The proof of \eqref{9e5.13} is now complete.
\end{proof}

 \section{Moments of the immigration density}\label{s3}
  
It remains to prove Lemma \ref{0l4.1}. Throughout the rest of the paper, we only consider $X_0\in M_c(\R)$ and $R>0$ as in \eqref{e1.04} following the assumption of Lemma \ref{0l4.1}. Recall from \eqref{4e2.4} that $ \dot{A}_t^R= V_t^R(f_R) $ where $V_t^R$ is as in \eqref{5e3} and
 \begin{align}\label{be3.41}
 f_R(x)= \int_{B_R^c}    \frac{c_\alpha}{|y-x|^{1+\alpha}} dy&=\frac{c_\alpha}{\alpha}\Big(\frac{1}{(R-x)^\alpha}+\frac{1}{(R+x)^\alpha}\Big), \quad \forall x\in B_R.
 \end{align}
 It follows that
  \begin{align}\label{4e6}
 f_R(x) \leq C(R-|x|)^{-\alpha}, \quad \forall x\in B_R.
 \end{align}
 Recall from \eqref{5e3} that $V^R$ is a superprocess whose spatial motion is the killed $\alpha$-stable process generated by $\Delta_\alpha^R$. Let $p_t^R(x, y)$ be the transition density of the killed $\alpha$-stable process and denote its the semigroup by $P_t^R$ such that for any function $\phi$,
  \begin{align*}
P_t^R \phi(x)=\int_{B_R} p_t^R(x, y) \phi(y) dy, \quad \forall x\in B_R.
  \end{align*}
 Theorem 1.1 of \cite{CKS10} states that for every $T>0$, the following holds for all $0<t\leq T$ and $x,y\in B_R$:
 \begin{align}\label{5e6}
p_t^R(x,y) \asymp \Big(1\wedge \frac{(R-|x|)^{\alpha/2}}{\sqrt{t}}\Big) \Big(1\wedge \frac{(R-|y|)^{\alpha/2}}{\sqrt{t}}\Big) \Big(t^{-1/\alpha}\wedge \frac{t}{|y-x|^{1+\alpha}} \Big).
  \end{align} 
  In the above, $f(t,x,y)\asymp g(t,x,y)$ means that $cg(t,x,y)\leq f(t,x,y)\leq Cg(t,x,y)$ for some constants $C>c>0$. 
 By using \eqref{5e6}, we will calculate the fourth moments for the difference of $V_t^R(f_R)$, thus proving Lemma \ref{0l4.1}.
 
   \begin{lemma} \label{l1.33}
   For any $\alpha \in (0, 2/3)$, there exist some constant $C>0$ depending only on $\alpha$ such that for any $0<\eps_0< 1\wedge (\frac{1}{\alpha}-\frac{3}{2})$, if $X_0\in M_c(\R)$ and $R>0$ satisfy \eqref{e1.04}, then
 \begin{align*} 
 \E\Big((V_{t+s}^R(f_R) -V_{t}^R(f_R) )^{4}\Big) \leq   Cs^{1+\eps_0} (X_0(1) \vee X_0(1)^4)  (R^{-\alpha} \vee R^{-8\alpha}), \quad \forall 0\leq t,s \leq 1.
\end{align*}
\end{lemma}

\begin{remark} \label{r4.2}
By replacing $f$ in \eqref{5e3} by $f_R$ (with some monotone arguments), we get
\begin{align*} 
 V_t^R(f_R)=X_0(f_R)+M^{V^R}_t(f_R)+\int_0^t   V_s^R(\Delta_\alpha^R f_R)  ds,
 \end{align*}
   where $M^{V^R}(f_R)$ is a continuous martingale with quadratic variation $t\mapsto \int_0^t  V^R_s(f_R^2) ds$. One can check by \eqref{be3.41} and \eqref{5e6} that the first moment of $\int_0^t V^R_s(f_R^2)ds$ is finite if and only if $\alpha\in (0, 2/3)$, thus justifying our condition for Lemma \ref{l1.33}. The moment calculations below also show that to obtain a continuous version of $(V_t^R(f_R), t\geq 0)$, one would require that $\alpha\in (0, 2/3)$. 
   
Note that the key step in our argument to prove that $W^R(1)$ as in \eqref{5e2.4} has lots of zeros is by comparing $W^R(1)$ with a square Bessel process. This in turn requires proving the continuity of $s\mapsto \dot{A}_s^R$. However, the continuity of $s\mapsto \dot{A}_s^R$ may not be required for $W^R(1)$ to still have lots of zeros, thus there is hope that our two main theorems could be extended to $\alpha\in (0, 2)$.
 \end{remark}

 Throughout the rest of the paper, we fix $0<\alpha<2/3$. In what follows, we write $\E_\mu$ for the law of $(V_t^R, t\geq 0)$ with $V_0^R=\mu$. Set 
 \begin{align*}
 \la \mu, f\ra=\mu(f)=\int f(x) \mu(dx).
  \end{align*}
  For any function $\phi\geq 0$ and $s>0$, we define 
\begin{align}\label{ae2.5}
 &v_{1}^\phi(s)=P_s^R \phi,\quad v_{n}^\phi(s)=\sum_{k=1}^{n-1} \binom{n-1}{k} \int_0^s P_{s-r}^R (v_{k}^\phi(u) v_{n-k}^\phi(u))  du, \ n\geq 2.
 \end{align}
 The following result on the moments of superprocesses is an easy consequence of Lemma 2.2 of \cite{KS88}.
 
\begin{lemma} \label{l4.0}
 For any $s>0$ and $\phi\geq 0$ such that $\|P_s^R \phi\|_\infty<\infty$, we have
 
 \no (i)
 \begin{align*}
 \E_\mu \la V_s^R, \phi \ra= \la \mu, v_{1}^\phi(s) \ra, \quad \quad\quad \E_\mu \la V_s^R, \phi \ra^2= \la \mu, v_{2}^\phi(s) \ra+\la \mu, v_{1}^\phi(s) \ra^2.
 \end{align*}
 (ii)
  \begin{align*}
 \E_\mu \la V_s^R, \phi \ra^3= \la \mu, v_{3}^\phi(s) \ra+3\la \mu, v_{2}^\phi(s) \ra\la \mu, v_{1}^\phi(s) \ra+\la \mu, v_{1}^\phi(s) \ra^3.
 \end{align*}
 (iii)
  \begin{align*}
 \E_\mu \la V_s^R, \phi \ra^4=  \la \mu, v_{4}^\phi(s) \ra&+4\la \mu, v_{3}^\phi(s) \ra\la \mu, v_{1}^\phi(s) \ra+3 \la \mu, v_{2}^\phi(s) \ra^2\\
 &+6\la \mu, v_{2}^\phi(s) \ra\la \mu, v_{1}^\phi(s) \ra^2+\la \mu, v_{1}^\phi(s) \ra^4
 \end{align*}
 (iv)
   \begin{align*}
 &\E_\mu (\la V_s^R,\phi \ra-\la \mu, P_s^R\phi \ra)^4= \la \mu, v_{4}^\phi(s) \ra+3 \la \mu, v_{2}^\phi(s) \ra^2.
 \end{align*}
  \end{lemma}
  \begin{proof}
  By Lemma 2.2 of \cite{KS88}, (i), (ii) and (iii) follow immediately by the monotone convergence theorem. (iv) is then an easy consequence by $v_{1}^\phi(s)=P_s^R \phi$.
  \end{proof}
  
The following technical result will be proved later in Appendix \ref{s4.3}.

\begin{lemma} \label{al1.0}
There exists some constant $C=C(\alpha)>0$ such that for any $0<\gamma< \frac{1}{\alpha}+\frac{1}{2}$, $0< s<1$ and $-R< y< R$, 
\begin{align*}
  \int_{B_R} p_{s}^R(y,x) \frac{1}{(R-|x|)^{\gamma\alpha}} dx  \leq  C\frac{1}{(R-|y|+s^{1/\alpha})^{\gamma\alpha}}. 
\end{align*}
\end{lemma}
Using similar calculus computations to the above, we obtain the following lemma, whose proof is deferred to Appendix \ref{s4.4}.
 \begin{lemma}\label{al1.1}
There exists some constant $C=C(\alpha)>0$ such that the following holds for all $0< s<1$ and $-R<x<R$:
\begin{align*}
  |P_{s}^R f_R(x)-f_R(x)| \leq &C\Big( (R-|x|)^{-\alpha} \wedge \frac{s}{(R-|x|)^{2\alpha}} \Big).
 \end{align*}
\end{lemma}
Lemma \ref{al1.0} ensures that we may apply Lemma \ref{l4.0} with $\phi=f_R$ as in \eqref{be3.41}.
To ease notation, when $\phi=f_R$, we write
  \begin{align}\label{be9.60}
v_n(s)=v_n^{f_R}(s)
 \end{align}
  for each $n\geq 1$ and any $s>0$.    Assuming the above two lemmas, we will prove Lemma \ref{l1.33} for $t=0$ and $t>0$ separately below. Fix $0<\eps_0< 1\wedge (\frac{1}{\alpha}-\frac{3}{2})$ throughout the rest of the section.

 \subsection{Proof of Lemma \ref{l1.33} for $t=0$}

Fix $s\in (0, 1)$. We first use $(a+b)^4\leq 2^4 a^4+2^4 b^4, \forall a,b\in \R$ to get
  \begin{align}\label{ce5.60}
&\E \Big(\Big(V_{s}^R(f_R)-V_0^R(f_R)\Big)^{4}\Big)\nn\\
& \leq 2^4 \E \Big(\Big(V_{s}^R(f_R)-\la V_0^R, P_s^R f_R\ra\Big)^{4}\Big)+2^4\Big(\la V_0^R, P_s^R f_R\ra-\la V_0^R, f_R\ra\Big)^{4}\nn\\
&\leq  C\la V_0^R, v_{4}(s) \ra+C \la V_0^R, v_{2}(s) \ra^2+C\la V_0^R, |P_s^R f_R-f_R|\ra^{4},
 \end{align}
 where the last inequality follows from Lemma \ref{l4.0} (iv) and \eqref{be9.60}. Using Lemma \ref{al1.1}, we obtain (recall $S(V_0^R)=S(X_0) \subseteq B(0,R/2)$)
   \begin{align}\label{be6.2}
\la V_0^R, |P_s^R f_R-f_R|\ra \leq \int_{|x|<R/2} Cs(R-|x|)^{-2\alpha} V_0^R(dx)\leq   CX_0(1) sR^{-2\alpha}.
   \end{align}
   It suffices to get bounds for $\la V_0^R, v_{4}(s) \ra$ and $\la V_0^R, v_{2}(s) \ra$. 
 
 Define
 \begin{align}\label{e5.60}
 \delta_0=\frac{1+\eps_0}{4} \in (0,1/2).
 \end{align}  
 One can easily check 
  \begin{align}\label{e5.61}
 1<4\delta_0= 1+\eps_0<1/\alpha-1/2<1/\alpha.
 \end{align}

 {\it Step 1.}  Recall $v_1(s,y)$ from \eqref{ae2.5} with $\phi=f_R$.  Apply \eqref{4e6} and Lemma \ref{al1.0} with $\gamma=1<\frac{1}{\alpha}+\frac{1}{2}$ to see that
  \begin{align}\label{ae2.7}
 v_1 (s,y)\leq C\int_{B_R} p_s^R(y,x) \frac{1}{(R-|x|)^\alpha} dx &\leq C  \frac{1}{(R-|y|+s^{1/\alpha})^{\alpha}}, \ \forall y\in B_R.
   \end{align}

 {\it Step 2.}  Use \eqref{ae2.5} with $n=2$ to get for all $y\in B_R$,
  \begin{align}\label{e2.17}
 v_2(s,y)&= \int_0^s du \int p_{s-u}^R(y,x)  v_1(u,x)^2 dx\nn\\
& \leq C   \int_0^s du \int p_{s-u}^R(y,x)  (R-|x|+u^{1/\alpha})^{-2\alpha} dx,
\end{align}
where we have used \eqref{ae2.7} in the inequality.

  The lemma below uses similar calculations to that of Lemmas \ref{al1.0} and \ref{al1.1}, so we defer the proof to Appendix \ref{4s2.1}. 

\begin{lemma} \label{l0.1}
There exists some constant $C=C(\alpha)>0$ such that for any $0< t<1$, $0\leq \gamma<1/\alpha-1/2$, $0\leq \rho\leq 1\wedge (1/\alpha-\gamma)$ and $-R< y< R$,
\begin{align*}
 &\int_0^t du \int p_{t-u}^R(y,x) \frac{1}{(R-|x|+u^{1/\alpha})^{(2+\gamma)\alpha}} dx \leq C \frac{t^{\rho}}{(R-|y|+t^{1/\alpha})^{(1+\gamma+\rho)\alpha}}. 
\end{align*}
\end{lemma}

Apply Lemma \ref{l0.1}  with $t=s,$ $\gamma=0$ and $\rho=2\delta_0\leq 1\wedge 1/\alpha$ to get that \eqref{e2.17} becomes
 \begin{align}\label{ae2.11}
v_2(s,y)&\leq    Cs^{2\delta_0} (R-|y|+s^{1/\alpha})^{-(1+2\delta_0)\alpha}, \ \forall y\in B_R.
\end{align}

{\it Step 3.}  Turning to $v_3(s,y)$ for $y\in B_R$, we use \eqref{ae2.7}, \eqref{ae2.11} and $u\leq s$ to get 
  \begin{align} \label{ae2.12}
 v_3(s,y)=&3\int_0^s du \int p_{s-u}^R(y,x)  v_1(u,x)  v_2(u,x) dx\nn\\
 \leq& Cs^{2\delta_0}  \int_0^s du \int p_{s-u}^R(y,x)  \frac{1}{(R-|x|+u^{1/\alpha})^{(2+ 2\delta_0)\alpha}} dx\nn\\
 \leq&    Cs^{3\delta_0} (R-|y|+s^{1/\alpha})^{-(1+3\delta_0)\alpha}.
   \end{align}
   In the last inequality, we have used Lemma \ref{l0.1} with $\gamma=2\delta_0<1/\alpha-1/2$ and \break  $\rho=\delta_0<1/\alpha-\gamma$.\\
   
 {\it Step 4.}    For $v_4(s,y)$ with $y\in B_R$, we get  
  \begin{align*}
 v_4(s,y)=&4\int_0^s du \int p_{s-u}^R(y,x)  v_1(u,x)  v_3(u,x) dx + 3\int_0^s du \int p_{s-u}^R(y,x)  v_2(u,x)^2  dx\nn\\
 \leq& C s^{3\delta_0} \int_0^s du \int p_{s-u}^R(y,x)  \frac{1}{(R-|x|+u^{1/\alpha})^{(2+ 3\delta_0)\alpha}} dx\nn\\
&\quad \quad + C s^{4\delta_0} \int_0^s du \int p_{s-u}^R(y,x)  \frac{1}{(R-|x|+u^{1/\alpha})^{(2+ 4\delta_0)\alpha}} dx.
   \end{align*}
The inequality above uses \eqref{ae2.7}, \eqref{ae2.11}, \eqref{ae2.12} and $u\leq s$. Apply Lemma \ref{l0.1} with $\gamma=3\delta_0<1/\alpha-1/2$ and $\rho=\delta_0<1/\alpha-\gamma$ to see that the first term above is bounded by $Cs^{4\delta_0} (R-|y|+s^{1/\alpha})^{-(1+4\delta_0)\alpha}.$ Next, apply Lemma \ref{l0.1} with $\gamma=4\delta_0<1/\alpha-1/2$ and $\rho=0$ to see that the second term above is also bounded by $Cs^{4\delta_0} (R-|y|+s^{1/\alpha})^{-(1+4\delta_0)\alpha}.$ Now we may conclude
  \begin{align}\label{ae2.13}
 v_4(s,y) \leq&   Cs^{4\delta_0} (R-|y|+s^{1/\alpha})^{-(1+4\delta_0)\alpha}, \ \forall y\in B_R.
   \end{align}   
 
  {\it Step 5.}  Using \eqref{ae2.11} and \eqref{ae2.13}, we obtain that (recall $S(V_0^R)=S(X_0) \subseteq B(0,R/2)$)
   \begin{align*} 
\la V_0^R, v_{2}(s) \ra \leq CX_0(1) s^{2\delta_0} R^{-(1+2\delta_0)\alpha}, \ \text{ and } \ \la V_0^R, v_{4}(s) \ra \leq CX_0(1) s^{4\delta_0} R^{-(1+4\delta_0)\alpha}.
   \end{align*}   
 Returning to \eqref{ce5.60}, we use the above and \eqref{be6.2} to conclude
   \begin{align*} 
\E ((V_{s}^R(f_R)-V_0^R(f_R))^{4}) \leq & CX_0(1)s^{4\delta_0} R^{-(1+4\delta_0)\alpha}+CX_0(1)^2s^{4\delta_0} R^{-(2+4\delta_0)\alpha}+ CX_0(1)^4s^4 R^{-8\alpha}\nn\\
\leq &C(X_0(1)\vee X_0(1)^4) (R^{-\alpha} \vee R^{-8\alpha}) s^{4\delta_0}.
 \end{align*}
 The proof is complete by $4\delta_0=1+\eps_0$.

\subsection{Proof of Lemma \ref{l1.33} for $t>0$}
We turn to the proof of Lemma \ref{l1.33} for $0<t, s<1$. Note that
  \begin{align}\label{be0.2}
\E\Big[(V_{t+s}^R(f_R)-V_t^R(f_R))^4\Big]&\leq 2^4 \E\Big[(V_{t+s}^R(f_R)-V_t^R(P_s^R f_R))^4\Big]\nn\\
&+2^4 \E\Big[(V_t^R(P_s^R f_R)-V_t^R(f_R))^4\Big]:=2^4 I_1+2^4 I_2.
\end{align}   
It suffices to bound $I_1$ and $I_2$.\\

{\bf Bounds for $I_1$}. By using the Markov property and Lemma \ref{l4.0} (iv), we get
  \begin{align}\label{be1.0}
I_1=& \E_{V_0^R}\Big[\E_{V_t^R}\Big[(V_{s}^R(f_R)-V_0^R(P_s^R f_R))^4\Big] \Big]\nn\\
=&\E_{V_0^R}\Big[\la V_t^R, v_{4}(s) \ra+ 3\la V_t^R, v_{2}(s) \ra^2\Big]\nn\\
= &\la V_0^R, v_{1}^{v_{4}(s)}(t) \ra+ 3\la V_0^R, v_{2}^{v_{2}(s)}(t) \ra+3\la V_0^R, v_{1}^{v_{2}(s)}(t) \ra^2,
\end{align}   
where the last equality follows from Lemma \ref{l4.0} (i) with $\phi=v_{4}(s)$ and $\phi=v_{2}(s)$. By \eqref{ae2.13}, we get for all $x\in B_R$,
  \begin{align*}
v_{1}^{v_{4}(s)}(t,x)=\int p_t^R(x,y) v_4(s,y) dy \leq&   Cs^{4\delta_0} \int p_t^R(x,y) (R-|y|+s^{1/\alpha})^{-(1+4\delta_0)\alpha} dy\nn\\
\leq &Cs^{4\delta_0} (R-|x|+t^{1/\alpha})^{-(1+4\delta_0)\alpha}.
   \end{align*}   
The last inequality above uses Lemma \ref{al1.0} with $1+4\delta_0<1/\alpha+1/2$ (recall \eqref{e5.61}). Hence by \eqref{e1.04},
  \begin{align}\label{be1.1}
   \la V_0^R, v_{1}^{v_{4}(s)}(t) \ra \leq  CX_0(1)  s^{4\delta_0} R^{-(1+4\delta_0)\alpha}.
   \end{align} 
 Similarly,  by \eqref{ae2.11}, we get for all $x\in B_R$,
    \begin{align}\label{be8.11}
v_{1}^{v_{2}(s)}(t,x) =\int p_t^R(x,y) v_2(s,y) dy \leq&   Cs^{2\delta_0} \int p_t^R(x,y) (R-|y|+s^{1/\alpha})^{-(1+2\delta_0)\alpha}dy\nn\\
\leq &Cs^{2\delta_0} (R-|x|+t^{1/\alpha})^{-(1+2\delta_0)\alpha},
   \end{align}  
      where the last inequality uses Lemma \ref{al1.0}. Next, by \eqref{ae2.5} with $\phi=v_2(s)$ and $n=2$, we get for all $x\in B_R$,
\begin{align}\label{be8.12}
v_{2}^{v_{2}(s)}(t,x) =& \int_0^t du \int p_{t-r}^R (x,y) v_{1}^{v_{2}(s)}(u,y)^2   dy\nn\\
\leq & Cs^{4\delta_0} \int_0^t du \int p_{t-u}^R(x,y) (R-|x|+u^{1/\alpha})^{-(2+4\delta_0)\alpha}  dy\nn\\
\leq &Cs^{4\delta_0}  (R-|x|+t^{1/\alpha})^{-(1+4\delta_0)\alpha},
   \end{align} 
   where the first inequality uses \eqref{be8.11}. The last inequality above follows by Lemma \ref{l0.1} with $\gamma=4\delta_0$ and $\rho=0$. Combine \eqref{e1.04}, \eqref{be8.11} and \eqref{be8.12} to get
     \begin{align}\label{be1.2}
  \la V_0^R, v_{1}^{v_{2}(s)}(t) \ra  &\leq  CX_0(1) s^{2\delta_0} R^{-(1+2\delta_0)\alpha}, \nn\\
   & \text{ and }  \la V_0^R, v_{2}^{v_{2}(s)}(t) \ra \leq  CX_0(1)s^{4\delta_0} R^{-(1+4\delta_0)\alpha}.
   \end{align} 
Returning to  \eqref{be1.0}, we apply \eqref{be1.1} and \eqref{be1.2}  to see that
  \begin{align}\label{be0.3}
I_1 &\leq CX_0(1) s^{4\delta_0} R^{-(1+4\delta_0)\alpha}+CX_0(1)s^{4\delta_0} R^{-(1+4\delta_0)\alpha}+ CX_0(1)^2 s^{4\delta_0} R^{-(2+4\delta_0)\alpha} \nn\\
&\leq  Cs^{4\delta_0} (X_0(1) \vee X_0(1)^4) (R^{-\alpha} \vee R^{-8\alpha}).\\\nn
\end{align}   

{\bf Bounds for $I_2$}. Turning to $I_2$, recall from \eqref{e5.60} that $\delta_0 \in (0,1)$. Use Lemma \ref{al1.1} and that $a\wedge b \leq a^{1-\delta_0} b^{\delta_0}, \forall a,b>0$ to see that 
 \begin{align*} 
 |P_s^R f_R(x)-f_R(x)| \leq  Cs^{\delta_0} (R-|x|)^{-(1+\delta_0)\alpha}, \quad \forall x\in B_R.
\end{align*}  
Define 
 \begin{align*} 
F_R(x)= (R-|x|)^{-(1+\delta_0)\alpha}, \quad \forall x\in B_R.
\end{align*}   
It follows that
  \begin{align}\label{ce0.1}
I_2=\E\Big[\Big(V_t^R(P_s^R f_R)-V_t^R(f_R)\Big)^4\Big]\leq Cs^{4\delta_0}\E(V_t^R(F_R)^4).
\end{align} 
Apply Lemma \ref{l4.0} (iii) with $\phi=F_R$ to get
  \begin{align}\label{be0.1}
\E(V_t^R(F_R)^4)\leq  &\la V_0^R, v_{4}^{F_R}(t) \ra+4\la V_0^R, v_{3}^{F_R}(t) \ra\la V_0^R, v_{1}^{F_R}(t) \ra\\
&+3 \la V_0^R, v_{2}^{F_R}(t) \ra^2+6\la V_0^R, v_{2}^{F_R}(t) \ra\la V_0^R, v_{1}^{F_R}(t) \ra^2+\la V_0^R, v_{1}^{F_R}(t) \ra^4.\nn
\end{align}   
It remains to get bounds for $v_n^{F_R}(t)$ for $1\leq n\leq 4$.\\

  {\it Step 1.} Recall $v_1^{F_R}(t)$ from \eqref{ae2.5} to get for all $-R<x<R$,
\begin{align} \label{be1.5}
 v_1^{F_R}(t,x)=\int p_t^R(x,y) F_R(y) dy&= \int p_t^R(x,y) \frac{1}{(R-|y|)^{(1+\delta_0)\alpha}}dy\nn\\
 &\leq \frac{C}{(R-|x|+t^{1/\alpha})^{(1+\delta_0)\alpha}},
\end{align}  
where the last inequality uses Lemma \ref{al1.0} with $\gamma=1+\delta_0<1/\alpha+1/2$. \\

  {\it Step 2.} Next, for all $-R<x<R$, we use \eqref{be1.5} to get
\begin{align} \label{be1.7}
 v_2^{F_R}(t,x) &=\int_0^t du\int p_{t-u}^R(x,y)   v_1^{F_R}(u,y)^2 dy\nn\\
 &\leq  C\int_0^t du\int p_{t-u}^R(x,y)   \frac{1}{(R-|y|+u^{1/\alpha})^{(2+2\delta_0)\alpha}} dy\nn\\
 &\leq C   \frac{1}{(R-|x|+t^{1/\alpha})^{(1+2\delta_0)\alpha}},
\end{align}  
where the last inequality follows by Lemma \ref{l0.1} with $\gamma=2\delta_0$ and $\rho=0$. \\

 {\it Step 3.} Turning to $v_3^{F_R}(t)$, we get  for all $-R<x<R$,
   \begin{align} \label{be1.9}
 v_3^{F_R}(t,x) &=3\int_0^t du\int p_{t-u}^R(x,y)   v_1^{F_R}(u,y)v_{2}^{F_R}(u,y) dy\nn\\
 &\leq  C\int_0^t du\int p_{t-u}^R(x,y)   \frac{1}{(R-|y|+u^{1/\alpha})^{(2+3\delta_0)\alpha}} dy\nn\\
 &\leq C   \frac{1}{(R-|x|+u^{1/\alpha})^{(1+3\delta_0)\alpha}},
\end{align}  
The second line above uses \eqref{be1.5} and \eqref{be1.7}. The third line follows by Lemma \ref{l0.1} with $\gamma=3\delta_0$ and $\rho=0$. \\

{\it Step 4.} To bound $v_4^{F_R}(t)$, we use \eqref{ae2.5} to see that  for all $-R<x<R$,
\begin{align} \label{be1.11}
 v_4^{F_R}(t,x) &=4\int_0^t du\int p_{t-u}^R(x,y)   v_1^{F_R}(u,y)v_{3}^{F_R}(u,y) dy\nn\\
&\quad+ 3\int_0^t du\int p_{t-u}^R(x,y)   v_{2}^{F_R}(u,y)^2 dy\nn\\
 &\leq  C\int_0^t du\int p_{t-u}^R(x,y)   \frac{1}{(R-|y|+u^{1/\alpha})^{(2+4\delta_0)\alpha}} dy\nn\\
 &\leq C   \frac{1}{(R-|x|+u^{1/\alpha})^{(1+4\delta_0)\alpha}},
\end{align}  
where the first inequality follows from \eqref{be1.5}, \eqref{be1.7}, and \eqref{be1.9}. The last inequality uses Lemma \ref{l0.1} with $\gamma=4\delta_0$ and $\rho=0$. \\

{\it Step 5.}  Recall $S(V_0^R)=S(X_0) \subseteq B(0,R/2)$ and use \eqref{be1.5} to see that
\begin{align*}
\la V_0^R, v_{1}^{F_R}(t) \ra \leq \int_{|x|<R/2} \frac{C}{(R-|x|+t^{1/\alpha})^{(1+\delta_0)\alpha}} V_0^R(dx) \leq CX_0(1)R^{-(1+\delta_0)\alpha},
\end{align*}   
Similarly, by \eqref{be1.7}, \eqref{be1.9} and \eqref{be1.11}, we get
\begin{align*}
&\la V_0^R, v_{2}^{F_R}(t) \ra \leq   CX_0(1) R^{-(1+2\delta_0)\alpha},\nn\\
&\la V_0^R, v_{3}^{F_R}(t) \ra \leq   CX_0(1) R^{-(1+3\delta_0)\alpha},\nn\\
&\la V_0^R, v_{4}^{F_R}(t) \ra \leq   CX_0(1) R^{-(1+4\delta_0)\alpha}.
\end{align*}   

Returning to \eqref{be0.1}, we conclude from the above that
  \begin{align*}
\E(V_t^R(F_R)^4)\leq  C(X_0(1)\vee X_0(1)^4) (R^{-\alpha}\vee R^{-8\alpha}).
\end{align*}   
Hence \eqref{ce0.1} becomes 
   \begin{align*}
I_2\leq Cs^{4\delta_0}(X_0(1)\vee X_0(1)^4) (R^{-\alpha}\vee R^{-8\alpha}).
\end{align*} 
 Together with \eqref{be0.2} and \eqref{be0.3}, we conclude that
  \begin{align*}
&\E\Big[(V_{t+s}^R(f_R)-V_t^R(f_R))^4\Big]\leq  Cs^{4\delta_0} (X_0(1) \vee X_0(1)^4)  (R^{-\alpha} \vee R^{-8\alpha}),
\end{align*}   
as required.

 \bibliographystyle{plain}
\def\cprime{$'$}

\appendix

 \section{Proof of Lemma \ref{l2.2}}\label{4a}
\begin{proof}[Proof of Lemma \ref{l2.2}]
Denote by $P^{(4\delta)}_x$ the law of the $4\delta$-dimensional square Bessel process $(Z_t, t\geq 0)$ starting from $x\geq 0$. The transition density of $Z$ under $P^{(4\delta)}_0$ is given by (see, e.g., Page 441 of \cite{RY94})
\[
q_t^{4\delta}(0,y)=\Gamma(2\delta)^{-1} (2t)^{-2\delta} y^{2\delta-1} \exp(-y/2t), \quad t>0, y>0.
\]
Set $T_0=\inf\{t> 0: Z_t=0\}$. 
By the Markov property, we get
\begin{align}\label{4e1.89}
I:=&P_0^{(4\delta)}(\text{Zeros}(Z) \cap (a, b)= \emptyset)=P_0^{(4\delta)}(P^{(4\delta)}_{{Z}_{a}} (T_0>b-a))\\
=&\int_0^\infty P^{(4\delta)}_{y} (T_0>b-a) \Gamma(2\delta)^{-1} (2a)^{-2\delta} y^{2\delta-1} \exp(-y/2a) dy.\nn
 \end{align}
 Set $\rho_t=\sqrt{Z_t}$ so that $\rho$ is a $4\delta$-dimensional Bessel process. Let $P^{4\delta}_x$ be the law of $\rho$ starting from $x\geq 0$. 
By slightly abusing the notation, we let $T_0=\inf\{t> 0: \rho_t=0\}$. One can easily check that 
\[
P^{(4\delta)}_{y} (T_0>b-a)=P^{4\delta}_{\sqrt{y}} (T_0>b-a).
\]
Then by using the density function of $T_0$ under $P^{4\delta}_x$ from Proposition 2.9 of Lawler \cite{Law18} to get (note the $a$ there is $2\delta-1/2$ by our notation)
\begin{align*}
& P^{4\delta}_{\sqrt{y}} (T_0>b-a)=\int_{b-a}^\infty c_\delta y^{1-2\delta} s^{2\delta-2} \exp{\{-y/2s\}} ds.
\end{align*}
Use the above in \eqref{4e1.89} to get
\begin{align}\label{4e1.39}
I&=\int_0^\infty  \Gamma(2\delta)^{-1} (2a)^{-2\delta} y^{2\delta-1} \exp(-y/2a) dy \int_{b-a}^\infty c_\delta y^{1-2\delta} s^{2\delta-2} \exp{\{-y/2s\}} ds\nn\\
&= C_\delta   a^{-2\delta}      \int_{b-a}^\infty   s^{2\delta-2}  ds \int_0^\infty     \exp(-y/2a) \exp{\{-y/2s\}}  dy,
 \end{align}
where the last equality uses Fubini's theorem. The integral of $y$ gives
\begin{align*}
\int_0^\infty     \exp(-y/2a) \exp{\{-y/2s\}}  dy =\frac{sa}{s+a}\leq a,
\end{align*}
and so
\begin{align*}
I &\leq   C_\delta   a^{1-2\delta}      \int_{b-a}^\infty  s^{2\delta-2}  ds \leq C_\delta   a^{1-2\delta}  (b-a)^{2\delta-1},
\end{align*}
as required. The proof is complete by recalling $I$ from \eqref{4e1.89}.
 \end{proof}

   \section{Proof of Lemma \ref{al1.0}}\label{s4.3}
This section is devoted to the proof of Lemma \ref{al1.0}. Set $0<\gamma< \frac{1}{\alpha}+\frac{1}{2}$.
Notice that for any $-R<x<R$, 
  \begin{align}\label{e0.78}
\frac{1}{(R-|x|)^{\gamma\alpha}} \leq \frac{1}{(R-x)^{\gamma\alpha}}+\frac{1}{(R+x)^{\gamma\alpha}}.
  \end{align}
It follows that for any $0< s<1$ and $-R< y< R$, 
  \begin{align}\label{e0.79}
\int p_{s}^R(y,x) \frac{1}{(R-|x|)^{\gamma\alpha}} dx &  \leq  \int p_{s}^R(y,x) \frac{1}{(R-x)^{\gamma\alpha}} dx   +\int p_{s}^R(y,x) \frac{1}{(R+x)^{\gamma\alpha}} dx\nn\\
&= \int p_{s}^R(y,x) \frac{1}{(R-x)^{\gamma\alpha}} dx   +\int p_{s}^R(-y,x) \frac{1}{(R-x)^{\gamma\alpha}} dx,
  \end{align}
 where the last equality follows by $p_{s}^R(y,-x)=p_{s}^R(-y,x)$.
 Hence it suffices to show
  \begin{align}\label{e0.01}
I:=\int p_{s}^R(y,x) \frac{1}{(R-x)^{\gamma\alpha}} dx  \leq C\frac{1}{(R-y+s^{1/\alpha})^{\gamma\alpha}}, \quad \forall -R<y<R.
  \end{align}

 \begin{proof}[Proof of Lemma \ref{al1.0}]

As noted above, the proof has been reduced to show \eqref{e0.01}. There are two cases for $y$: (1) $R-2s^{1/\alpha}<y<R$; (2) $-R<y<R-2s^{1/\alpha}$. We will prove that \eqref{e0.01} holds for both cases.\\

  {\bf Case 1.} We first consider $R-2s^{1/\alpha}<y<R$, in which case $s^{-\gamma}\leq 3^{\gamma \alpha}(R-y+s^{1/\alpha})^{-\gamma\alpha}$. Hence it suffices to show that
  \begin{align*}
I=\int  p_s^R(y,x) \frac{1}{(R-x)^{\gamma\alpha}} dx \leq Cs^{-\gamma}, \quad \forall R-2s^{1/\alpha}<y<R.
  \end{align*}

  \no Fix $R-2s^{1/\alpha}<y<R$. Separate the integral in $I$ into three integrals over three regions: (i) $R-s^{1/\alpha}<x<R$; (ii) $y-s^{1/\alpha}< x< R-s^{1/\alpha}$; (iii) $-R< x< y-s^{1/\alpha}$. Denote each corresponding integral by $I^{(i)}, I^{(ii)}, I^{(iii)}$ respectively so that $I=I^{(i)}+I^{(ii)}+I^{(iii)}$.  
  \\
 
   {\bf Part 1(i).}   
For the case of $R-s^{1/\alpha}<x<R$,  by \eqref{5e6} we get
 \begin{align*} 
p_s^R(y,x) \leq C   \frac{(R-x)^{\alpha/2}}{\sqrt{s}} s^{-1/\alpha}.
  \end{align*}
  It follows that
 \begin{align*}
I^{(i)}= \int_{R-s^{1/\alpha}}^R p_s^R(y,x) \frac{1}{(R-x)^{\gamma\alpha}} dx 
  &\leq C s^{-1/\alpha -1/2} \int_{R-s^{1/\alpha}}^R (R-x)^{-(\gamma-1/2)\alpha}   dx \\
 &\leq C  s^{-1/\alpha -1/2} (s^{1/\alpha})^{1-(\gamma-1/2)\alpha}=C  s^{-\gamma},
 \end{align*}
 where the second inequality follows by $\gamma<1/\alpha+1/2$.\\
 
  {\bf Part 1(ii).} Consider the case of $y-s^{1/\alpha}< x< R-s^{1/\alpha}$. By $R-x>s^{1/\alpha}$ and \eqref{5e6}, we get
 \begin{align*} 
\frac{1}{(R-x)^{\gamma\alpha}}\leq s^{-\gamma}  \quad \text{ and }\quad   p_s^R(y,x) \leq C   s^{-1/\alpha}.
  \end{align*}
  It follows that
 \begin{align*}
 I^{(ii)}=&\int_{y-s^{1/\alpha}}^{R-s^{1/\alpha}} p_s^R(y,x) \frac{1}{(R-x)^{\gamma\alpha}} dx\leq C    s^{-1/\alpha} s^{-\gamma} (R-y)\leq Cs^{-\gamma},
 \end{align*}
 where the last inequality follows by $R-y<2s^{1/\alpha}$.\\

  {\bf Part 1(iii).} If $-R\leq x\leq y-s^{1/\alpha}$, by \eqref{5e6} and $y<R$ we have
 \begin{align*} 
p_s^R(y,x) \leq C   \frac{s}{(y-x)^{1+\alpha}} \ \text{ and } \  \frac{1}{(R-x)^{\gamma\alpha}}\leq \frac{1}{(y-x)^{\gamma\alpha}}.
  \end{align*}
  It follows that
 \begin{align*}
I^{(iii)}&= \int_{-R}^{y-s^{1/\alpha}}  p_s^R(y,x) \frac{1}{(R-x)^{\gamma\alpha}} dx  \leq C\int_{-R}^{y-s^{1/\alpha}}      \frac{s}{(y-x)^{1+\alpha}} \frac{1}{(y-x)^{\gamma\alpha}}dx\\
&=   C  s \int_{s^{1/\alpha}}^{y+R}    \frac{1}{r^{1+(\gamma+1)\alpha}}  dr\leq Cs (s^{1/\alpha})^{-(\gamma+1)\alpha}=Cs^{-\gamma}.
 \end{align*}
 
The proof of \eqref{e0.01} for $R-2s^{1/\alpha}<y<R$ is now complete by the above three parts.\\

      {\bf Case 2.}   Turning to the case when $-R<y<R-2s^{1/\alpha}$, we get  
\begin{align*}
R-y>  \frac{R-y+s^{1/\alpha}}{2}, \text{ and so } (R-y)^{-\gamma\alpha}\leq 2^{\gamma \alpha}(R-y+s^{1/\alpha})^{-\gamma\alpha}.
\end{align*}     
Thus, it suffices to show that
  \begin{align*}
I=\int  p_s^R(y,x) \frac{1}{(R-x)^{\gamma\alpha}} dx \leq C(R-y)^{-\gamma \alpha}, \quad \forall -R<y<R-2s^{1/\alpha}.
  \end{align*}
Fix $-R<y<R-2s^{1/\alpha}$. Separate the integral in $I$ into four integrals over four regions: (i) $(R+y)/2<x<R$; (ii) $y+s^{1/\alpha}< x< (R+y)/2$; (iii) $y-s^{1/\alpha}< x< y+s^{1/\alpha}$; (iv) $-R<x< y-s^{1/\alpha}$.  Denote each corresponding integral by $I^{(i)}, I^{(ii)}, I^{(iii)}, I^{(iv)}$ respectively so that $I=I^{(i)}+I^{(ii)}+I^{(iii)}+I^{(iv)}$. \\
 
   {\bf Part 2(i).} If $(R+y)/2<x<R$, we get $|x-y|=x-y>  (R-y)/2$. Use this and \eqref{5e6} to see that
 \begin{align*} 
p_s^R(y,x) \leq C   \frac{(R-x)^{\alpha/2}}{\sqrt{s}} \frac{s}{|x-y|^{1+\alpha}}\leq Cs^{1/2} (R-x)^{\alpha/2} \frac{1}{(R-y)^{1+\alpha}}.
  \end{align*}
  It follows that
 \begin{align*}
 I^{(i)} =\int_{(R+y)/2}^R p_s^R(y,x)\frac{1}{(R-x)^{\gamma\alpha}} dx& \leq  \frac{Cs^{1/2}}{(R-y)^{1+\alpha}}    \int_{(R+y)/2}^R (R-x)^{-(\gamma-1/2)\alpha}   dx\\
&= \frac{Cs^{1/2}}{(R-y)^{1+\alpha}}    \int_0^{(R-y)/2}  r^{-(\gamma-1/2)\alpha}   dr.
 \end{align*}
Notice that $(\gamma-1/2)\alpha<1$ by $\gamma<1/\alpha+1/2$. So the above is at most
 \begin{align*}
I^{(i)}\leq &  C\frac{s^{1/2}}{(R-y)^{1+\alpha}}  (R-y)^{1-(\gamma-1/2)\alpha} \leq C\frac{1}{(R-y)^{\gamma\alpha}},
 \end{align*}
 where the last inequality uses $(R-y)>2s^{1/\alpha}$. \\
 
  {\bf Part 2(ii).}  Consider the case of $y+s^{1/\alpha}\leq x\leq (R+y)/2$, we get $R-x\geq (R-y)/2$ and so
 \begin{align*} 
 \frac{1}{(R-x)^{\gamma\alpha}}\leq C\frac{1}{(R-y)^{\gamma\alpha}}.
  \end{align*}
  Use the above and $p_s^R(y,x) \leq Cs/{|x-y|^{1+\alpha}}$ from \eqref{5e6} to see that
 \begin{align*}
 I^{(ii)}=\int_{y+s^{1/\alpha}}^{(R+y)/2} p_s^R(y,x) \frac{1}{(R-x)^{\gamma\alpha}} dx  &  \leq   Cs (R-y)^{-\gamma\alpha} \int_{y+s^{1/\alpha}}^{(R+y)/2}  \frac{1}{(x-y)^{1+\alpha}}  dx\\
 & \leq Cs (R-y)^{-\gamma\alpha} (s^{1/\alpha})^{-\alpha}=C(R-y)^{-\gamma\alpha}.
 \end{align*}

  {\bf Part 2(iii).}   If $y-s^{1/\alpha}\leq x\leq y+s^{1/\alpha}$, we get $R-x\geq R-y-s^{1/\alpha}\geq (R-y)/2$ and so
 \begin{align*} 
\frac{1}{(R-x)^{\gamma\alpha}}\leq C\frac{1}{(R-y)^{\gamma\alpha}}.
  \end{align*}
  By using the above and   $p_s^R(y,x) \leq C   s^{-1/\alpha}$ from \eqref{5e6}, we obtain
  \begin{align*}
 I^{(iii)}=&\int_{y-s^{1/\alpha}}^{y+s^{1/\alpha}} p_s^R(y,x) \frac{1}{(R-x)^{\gamma\alpha}} dx \leq \int_{y-s^{1/\alpha}}^{y+s^{1/\alpha}}  C s^{-1/\alpha}  \frac{1}{(R-y)^{\gamma\alpha}} dx\leq C\frac{1}{(R-y)^{\gamma\alpha}}.
 \end{align*}

   {\bf Part 2(iv).} Finally for the case of $-R\leq x\leq y-s^{1/\alpha}$, by \eqref{5e6} and $x<y$ we get  
 \begin{align*} 
p_s^R(y,x) \leq C   \frac{s}{(y-x)^{1+\alpha}} \text{ and }   \frac{1}{(R-x)^{\gamma\alpha}}\leq \frac{1}{(R-y)^{\gamma\alpha}}.
  \end{align*}
  It follows that
 \begin{align*} 
I^{(iv)}= \int_{-R}^{y-s^{1/\alpha}}  p_s^R(y,x)  \frac{1}{(R-x)^{\gamma\alpha}} dx  & \leq Cs\frac{1}{(R-y)^{\gamma\alpha}} \int_{-R}^{y-s^{1/\alpha}}      \frac{1}{(y-x)^{1+\alpha}} dx\\
& \leq Cs\frac{1}{(R-y)^{\gamma\alpha}} (s^{1/\alpha})^{-\alpha}= C \frac{1}{(R-y)^{\gamma\alpha}}.
 \end{align*}
 
The proof of \eqref{e0.01} for $-R<y<R-2s^{1/\alpha}$ is now complete given the four parts above. 
\end{proof}

  \section{Proof of Lemma \ref{al1.1}} \label{s4.4}
Recall $f_R$ from \eqref{be3.41}.
Define 
  \begin{align}
g_R(z)=\frac{1}{(R-z)^\alpha} \text{ and } h_R(z)=\frac{1}{(R+z)^\alpha}, \quad \forall -R<z<R.
 \end{align}
 Then for any $x\in B_R$, $f_R(x)=c_\alpha \alpha^{-1} (g_R(x)+h_R(x))$. It follows that for any $s>0$,
     \begin{align}\label{e0.03}
 P_{s}^R f_R(x)=c_\alpha \alpha^{-1}  P_{s}^R g_R(x)+c_\alpha \alpha^{-1}  P_{s}^R h_R(x), \quad \forall -R<x<R.
  \end{align}
  Since $h_R(y)=g_R(-y), \forall y\in B_R$, by symmetry we obtain
\[
P_{s}^R h_R(x)=\int p_{s}^R(x,y) h_R(y) dy=\int p_{s}^R(-x,y) g_R(y) dy=P_{s}^R g_R(-x), \ \ \forall x\in B_R.
\]
It follows that for any $s>0$ and $x\in B_R$,
\begin{align*}
  |P_{s}^R f_R(x)-f_R(x)|\leq&c_\alpha \alpha^{-1}   |P_{s}^R g_R(x)-g_R(x)|+c_\alpha \alpha^{-1}   |P_{s}^R h_R(x)-h_R(x)|\\
  =&c_\alpha\alpha^{-1}   |P_{s}^R g_R(x)-g_R(x)|+c_\alpha \alpha^{-1}   |P_{s}^R g_R(-x)-g_R(-x)|.
\end{align*}
 Hence it suffices to show that for any $s>0$, 
 \begin{align}\label{e2.01}
 |P_{s}^R g_R(x)-g_R(x)|\leq C \Big( (R-|x|)^{-\alpha} \wedge \frac{s}{(R-|x|)^{2\alpha}} \Big), \ \forall x\in B_R.
  \end{align}

 To do this, we notice that for $x\in B_R$,
   \begin{align*}
 P_{s}^R g_R(x)-g_R(x)=\int_{-R}^R p_s^R(x,z)  (g_R(z)-g_R(x))dz -g_R(x) \Big(1-\int_{-R}^R p_s^R(x,z)    dz\Big).
  \end{align*}
   Define
\begin{align}
I_s(x):=\int_{-R}^R p_s^R(x,z) |g_R(z)-g_R(x)| dz, \quad \forall x\in B_R.
\end{align}
Then we get for each $x\in B_R$,
 \begin{align} \label{be4.1}
 |P_{s}^R g_R(x)-g_R(x)|\leq I_s(x) +   g_R(x) \Big(1-\int_{-R}^R p_s^R(x,z)    dz\Big).
  \end{align}
  Let $Y=\{Y_t, t\geq 0\}$ be a symmetric $\alpha$-stable process on $\R$ starting from $r\in \R$, the law of which is denoted by $P^r$. Then for the second term above, we have
 \begin{align} \label{be4.2}
&1-\int_{-R}^R p_s^R(x,z)    dz= P^x\Big(\sup_{0\leq r\leq s} |Y_r|>R\Big)\leq 2 P^0\Big(\sup_{0\leq r\leq s} Y_r>R-|x|\Big), \quad \forall x\in B_R.
\end{align}
    Proposition VIII.4 of \cite{Ber96}  gives us
 \begin{align} \label{be7.1}
 \lim_{r\to \infty } r^\alpha P^0\Big(\sup_{0\leq u\leq 1} Y_u \geq r\Big) =C
 \end{align}
  for some constant $C>0$.
  
Recall the scaling property of $Y$ that for any $t>0$, $(t^{-\frac{1}{\alpha}}Y_{ut}, u\geq 0)$ has the same distribution as $(Y_u, u\geq 0)$. One may get
 \begin{align} \label{be4.3}
   &P^0\Big(\sup_{0\leq r\leq s} Y_r>R-|x|\Big)=P^0\Big(\sup_{0\leq u\leq 1} s^{-1/\alpha}Y_{us}>s^{-1/\alpha}(R-|x|)\Big)\nn\\
   &=P^0\Big(\sup_{0\leq u\leq 1}  Y_{u}>s^{-1/\alpha}(R-|x|)\Big)\leq 1\wedge  \Big(C \frac{s}{(R-|x|)^\alpha}\Big), \quad \forall x\in B_R,
\end{align}
where the last inequality uses \eqref{be7.1}. Combine \eqref{be4.1}, \eqref{be4.2} and \eqref{be4.3} to conclude that
 \begin{align*} 
 |P_{s}^R g_R(x)-g_R(x)|\leq I_s(x) + C g_R(x) \Big(1\wedge   \frac{s}{(R-|x|)^\alpha}\Big), \quad \forall x\in B_R.
  \end{align*}
Since $g_R(x)\leq (R-|x|)^{-\alpha}$,  the proof of \eqref{e2.01} is now reduced to show that
 \begin{align} \label{be7.4}
I_s(x) \leq C \Big( (R-|x|)^{-\alpha} \wedge \frac{s}{(R-|x|)^{2\alpha}} \Big),   \quad \forall x\in B_R.
  \end{align}

  \begin{proof}[Proof of Lemma \ref{al1.1}]
Fix any $s\in (0,1)$. As noted above, it suffices to prove \eqref{be7.4}.
If $R-2s^{1/\alpha}\leq |x|<R$, then $(R-|x|)^{-\alpha} \leq  2^\alpha {s}{(R-|x|)^{-2\alpha}}$. So it suffices to show that
\[
I_s(x)=\int_{-R}^R p_s^R(x,z) |g_R(z)-g_R(x)| dz\leq C(R-|x|)^{-\alpha}, \quad\forall    R-2s^{1/\alpha}\leq |x|<R.
\]
 To see this, by Lemma \ref{al1.0} with $\gamma=1$, we get for $R-2s^{1/\alpha}\leq |x|<R$,
\[
\int_{-R}^R p_s^R(x,z) g_R(z) dz\leq \int_{-R}^R p_s^R(x,z)  \frac{1}{(R-|z|)^\alpha}  dz\leq C(R-|x|)^{-\alpha}.
\]
It follows that
\begin{align*}
I_s(x)& \leq \int_{-R}^R p_s^R(x,z)  g_R(z) dz+g_R(x) \int_{-R}^R p_s^R(x,z)   dz\\
&\leq \frac{C}{(R-|x|)^\alpha} +\frac{1}{(R-x)^\alpha}\leq C(R-|x|)^{-\alpha}, \quad\forall    R-2s^{1/\alpha}\leq |x|<R,
\end{align*}
 as required. 
 
 It remains to consider the case of $|x|<R-2s^{1/\alpha}$, in which  ${s}{(R-|x|)^{-2\alpha}} \leq (R-|x|)^{-\alpha}$. Thus it is sufficient to show that
\begin{align}\label{be7.9}
I_s(x)\leq Cs(R-|x|)^{-2\alpha}, \quad\forall  |x|<R-2s^{1/\alpha}.
  \end{align}

   For any $x, z\in (-R,R)$, we get
\begin{align}\label{4e1.00}
   |g_R(z)-g_R(x)|= \frac{|(R-x)^\alpha-(R-z)^\alpha|}{(R-z)^\alpha(R-x)^\alpha} =\alpha   {|z-x|} \frac{(R-\xi)^{\alpha-1}}{(R-z)^\alpha(R-x)^\alpha},
  \end{align}
  where the last inequality follows by the mean-value theorem and $\xi$ is between $x$ and $z$.  Recall that $\alpha<2/3<1$. If $z>x$, we get $R-z<R-\xi<R-x$, so
\begin{align}\label{4e1.01}
    &|g_R(z)-g_R(x)|\leq C|x-z|\frac{1}{(R-z)(R-x)^\alpha}.
  \end{align}
 If $z<x$, we get $R-x<R-\xi<R-z$, so
\begin{align}\label{4e1.02}
 &|g_R(z)-g_R(x)|\leq C|x-z|\frac{1}{(R-z)^\alpha(R-x)}.
\end{align}

To prove \eqref{be7.9} for $-R+2s^{1/\alpha}<x<R-2s^{1/\alpha}$, we separate the integral in $I_s(x)$ into five integrals over five regions:   (i) $R-s^{1/\alpha}<z<R$; 
   (ii) $(R+x)/2<z<R-s^{1/\alpha}$; 
  (iii) $x+s^{1/\alpha}<z<(R+x)/2$; 
  (iv) $x-s^{1/\alpha}<z<x+s^{1/\alpha}$;  
  (v) $-R<z<x-s^{1/\alpha}$. Denote each corresponding integral by $I^{(i)}, I^{(ii)}, I^{(iii)}, I^{(iv)}, I^{(v)}$ respectively so that $I_s(x)=I^{(i)}+ I^{(ii)}+I^{(iii)}+I^{(iv)}+I^{(v)}$.
\\
      
{\bf Part (i).} Consider the case of $R-s^{1/\alpha}<z<R$. By \eqref{5e6}, we get 
 \begin{align*} 
p_s^R(x,z) \leq C   \frac{(R-z)^{\alpha/2}}{\sqrt{s}} \frac{s}{|z-x|^{1+\alpha}}.
  \end{align*}
 It follows that
 \begin{align*}
 I^{(i)}=&\int_{R-s^{1/\alpha}}^R p_s^R(x,z) |g_R(z)-g_R(x)| dz \leq C  s^{1/2} \int_{R-s^{1/\alpha}}^R \frac{(R-z)^{\alpha/2}    }{|x-z|^{1+\alpha}}  |g_R(z)-g_R(x)| dz.
 \end{align*}
 Use \eqref{4e1.01} with $z>x$ to see that  the above is at most 
  \begin{align*}
  I^{(i)}&\leq  C  s^{1/2} \int_{R-s^{1/\alpha}}^R   \frac{(R-z)^{\alpha/2}  }{|x-z|^\alpha}\frac{1}{(R-z)(R-x)^\alpha} dz.
    \end{align*}
    Noticing that $|x-z|=z-x\geq R-s^{1/\alpha}-x\geq (R-x)/2$, we get
    \begin{align*}
  I^{(i)}& \leq  C  s^{1/2} (R-x)^{-2\alpha}\int_{R-s^{1/\alpha}}^R   \frac{1}{(R-z)^{1-\alpha/2} } dz\\
    &\leq  C  s^{1/2} (R-x)^{-2\alpha} (s^{1/\alpha})^{\alpha/2}=Cs (R-x)^{-2\alpha}.
    \end{align*}

  {\bf Part (ii).} For the case of $(R+x)/2<z<R-s^{1/\alpha}$, we  get  $|z-x|\geq  (R-x)/2$ and thus by \eqref{5e6},
 \begin{align*} 
p_s^R(x,z) \leq C \frac{s}{|z-x|^{1+\alpha}}\leq C \frac{s}{(R-x)^{1+\alpha}}.
  \end{align*}
Use the above to see that
   \begin{align*}
 I^{(ii)} &=\int_{(R+x)/2}^{R-s^{1/\alpha}} p_s^R(x,z) |g_R(z)-g_R(x)| dz  \leq C \frac{s}{(R-x)^{1+\alpha}} \int_{(R+x)/2}^{R-s^{1/\alpha}}      (g_R(z)+g_R(x)) dz.
  \end{align*}
 Noticing that 
    \begin{align*}
 \int_{(R+x)/2}^{R-s^{1/\alpha}}   g_R(z) dz= \int_{(R+x)/2}^{R-s^{1/\alpha}}   \frac{1}{(R-z)^\alpha} dz=\int_{s^{1/\alpha}}^{(R-x)/2}    r^{-\alpha} dr\leq C(R-x)^{1-\alpha},
  \end{align*}
 and 
 \begin{align*}
 \int_{(R+x)/2}^{R-s^{1/\alpha}}  g_R(x) dz=\frac{1}{(R-x)^\alpha}  \Big(\frac{R-x}{2}-s^{1/\alpha}\Big)\leq (R-x)^{1-\alpha},
   \end{align*}
   we obtain
   
  \begin{align*}
 I^{(ii)} & \leq C \frac{s}{(R-x)^{1+\alpha}} \cdot  C(R-x)^{1-\alpha} \leq Cs (R-x)^{-2\alpha}.
 \end{align*}
 
   {\bf Part (iii).} Turning to $x+s^{1/\alpha}<z<(R+x)/2$, we apply \eqref{5e6} to get  
 \begin{align*} 
p_s^R(x,z) \leq C \frac{s}{|z-x|^{1+\alpha}}.
  \end{align*}
Use \eqref{4e1.01} with $z>x$ and the above to see  
 \begin{align*}
 I^{(iii)}=&\int_{x+s^{1/\alpha}}^{(R+x)/2} p_s^R(x,z) |g_R(z)-g_R(x)| dz \leq C  \int_{x+s^{1/\alpha}}^{(R+x)/2} \frac{s}{|z-x|^{\alpha}}  \frac{1}{(R-z)(R-x)^\alpha} dz.
 \end{align*}
By using $R-z>(R-x)/2$, we obtain 
         \begin{align*}
  I^{(iii)} &\leq  C\frac{s}{(R-x)^{1+\alpha}}   \int_{x+s^{1/\alpha}}^{(R+x)/2} \frac{1}{(z-x)^\alpha  } dz\\
   & \leq C   \frac{s}{(R-x)^{1+\alpha}} \cdot C(R-x)^{1-\alpha} \leq C   \frac{s}{(R-x)^{2\alpha}}.
    \end{align*}

   {\bf Part (iv).} For $z\in (x-s^{1/\alpha}, x+s^{1/\alpha})$, by \eqref{5e6} we get  
 \begin{align} \label{e1.37}
p_s^R(x,z) \leq C  s^{-1/\alpha}.
  \end{align}
  Next, we use \eqref{4e1.00} to see that
 \begin{align} \label{be1.38}
  |g_R(z)-g_R(x)| =\alpha  |x-z|\frac{|R-\xi|^{\alpha-1}}{(R-z)^\alpha(R-x)^\alpha}  \leq   \frac{s^{1/\alpha}}{(R-z)^\alpha(R-x)^\alpha (R-\xi)^{1-\alpha}},
 \end{align}
 where inequality follows since $|x-z|\leq s^{1/\alpha}$ and $\alpha<2/3$. Recall that $\xi$ is between $z$ and $x$. Hence 
  \begin{align*}
  R-z\geq R-x-s^{1/\alpha}\geq (R-x)/2  \ \text{ and } \ R-\xi \geq R-x-s^{1/\alpha}\geq (R-x)/2.
   \end{align*}
Apply the above in \eqref{be1.38} to get
 \begin{align} \label{e1.38}
  |g_R(z)-g_R(x)| \leq Cs^{1/\alpha}  {(R-x)^{-1-\alpha}}.
 \end{align}
Combining \eqref{e1.37} and \eqref{e1.38}, we conclude that
 \begin{align*}
 I^{(iv)}&=\int_{x-s^{1/\alpha}}^{x+s^{1/\alpha}} p_s^R(x,z) |g_R(z)-g_R(x)| dz \leq  \int_{x-s^{1/\alpha}}^{x+s^{1/\alpha}}   C  s^{-1/\alpha}  \cdot Cs^{1/\alpha}  {(R-x)^{-1-\alpha}} dz\\
 &\leq  C s^{1/\alpha}(R-x)^{-\alpha-1} \leq Cs(R-x)^{-2\alpha},
    \end{align*}
    where the last inequality follows by $\alpha<1$ and $R-x>2s^{1/\alpha}$.\\

  {\bf Part (v).}   Finally for the case of $-R\leq z\leq x-s^{1/\alpha}$, we use \eqref{5e6} to obtain  
 \begin{align*} 
p_s^R(x,z) \leq C    \frac{s}{(x-z)^{1+\alpha}},
  \end{align*}
and hence
 \begin{align} \label{be7.10}
  I^{(v)}&=\int_{-R}^{x-s^{1/\alpha}} p_s^R(x,z) |g_R(z)-g_R(x)| dz \nn\\
  &\leq Cs \int_{-R}^{x-s^{1/\alpha}} \frac{1}{(x-z)^{1+\alpha}} |{(R-z)^{-\alpha}} -{(R-x)^{-\alpha}} | dz.
  \end{align}
  
  Recall that $-R+2s^{1/\alpha}<x<R-2s^{1/\alpha}$. We first consider $-R+2s^{1/\alpha}< x\leq 0$. Use \eqref{4e1.02} for $z<x$ to see that 
\begin{align*}
  I^{(v)} & \leq C   s\int_{-R}^{x-s^{1/\alpha}}    \frac{1}{(x-z)^{\alpha}} \frac{1}{(R-z)^\alpha(R-x)} dz\\
  &\leq C\frac{s}{(R-x)^{1+\alpha}}\int_{-R}^{x-s^{1/\alpha}}    \frac{1}{(x-z)^{\alpha}} dz\\
  &\leq C\frac{s}{(R-x)^{1+\alpha}} \cdot C (R+x)^{1-\alpha} \leq Cs(R-|x|)^{-2\alpha},
  \end{align*}
  where the second inequality uses $R-z>R-x$. The last inequality follows by $R+x=R-|x|$ and $R-x\geq R-|x|$ for $x\leq 0$.\\

Next, we consider $0<x<R-2s^{1/\alpha}$. Recall \eqref{be7.10} to see that
  \begin{align*}
  I^{(v)} \leq &  Cs \int_{-R}^{-R+2x}  \frac{1}{(x-z)^{1+\alpha}} |{(R-z)^{-\alpha}} -{(R-x)^{-\alpha}}| dz \\
&+Cs\int_{-R+2x}^{x-s^{1/\alpha}}  \frac{1}{(x-z)^{1+\alpha}} |{(R-z)^{-\alpha}} -{(R-x)^{-\alpha}}| dz:=J_1+J_2.
     \end{align*}
 For the first term $J_1$, since $z<-R+2x$, we get $${(R-z)^{-\alpha}}\leq (R-(-R+2x))^{-\alpha}= (2R-2x)^{-\alpha},$$  and hence $|{(R-z)^{-\alpha}} -{(R-x)^{-\alpha}}| \leq C(R-x)^{-\alpha}$. It follows that
 \begin{align} \label{be7.11}
J_1&\leq   C  s (R-x)^{-\alpha} \int_{-R}^{-R+2x}   \frac{1}{(x-z)^{1+\alpha}}  dz  \leq  Cs (R-x)^{-2\alpha}.
\end{align} 
 For the second term $J_2$, we use \eqref{4e1.02} with $z<x$ to get
\begin{align*}
J_2 & \leq C   s\int_{-R+2x}^{x-s^{1/\alpha}}    \frac{1}{(x-z)^{\alpha}} \frac{1}{(R-z)^\alpha(R-x)} dz.
  \end{align*}
Use $R-z\geq R-x$ to further obtain
 \begin{align} \label{be7.12}
J_2 & \leq    C  s (R-x)^{-1-\alpha} \int_{-R+2x}^{x-s^{1/\alpha}}     \frac{1}{(x-z)^{\alpha}}  dz\nn\\
&=    C  s (R-x)^{-1-\alpha}   \int_{s^{1/\alpha}}^{R-x}    \frac{1}{r^{ \alpha}}  dr\leq Cs (R-x)^{-2\alpha}.
\end{align}
Combine \eqref{be7.11} and \eqref{be7.12} to conclude $I^{(v)}\leq J_1+J_2\leq Cs (R-x)^{-2\alpha}$ for $x>0$. \\

The proof of \eqref{be7.9} is now complete given the above five parts.
 \end{proof}

    \section{Proof of Lemmas \ref{l0.1} } \label{4s2.1}

\begin{proof}[Proof of Lemma \ref{l0.1}]
Fix $0< t<1$ and $0\leq \gamma<1/\alpha-1/2$. For each $-R< y< R$, define
 \[
 I:= \int_0^t du \int p_{t-u}^R(y,x) \frac{1}{(R-x+u^{1/\alpha})^{(2+\gamma)\alpha}} dx.
 \]
 We claim that it suffices to show that
 \begin{align}\label{e1.56}
 I\leq C  \frac{t}{(R-y+t^{1/\alpha})^{ (2+\gamma) \alpha}}+C \frac{t^{1/\alpha-\gamma}}{(R-y+t^{1/\alpha})^{1+\alpha}}, \quad \forall y\in B_R.
 \end{align}
Assuming the above, for any $0\leq \rho \leq 1\wedge (1/\alpha-\gamma)$, if $R-y<2t^{1/\alpha}$, by \eqref{e1.56} we get
  \begin{align*} 
 I\leq C   \frac{t}{(t^{1/\alpha})^{ (2+\gamma) \alpha}}+C \frac{t^{1/\alpha-\gamma}}{(t^{1/\alpha})^{1+\alpha}}\leq Ct^{-1-\gamma} \leq C\frac{t^\rho}{(R-y+t^{1/\alpha})^{ (1+\gamma+\rho) \alpha}},
 \end{align*}
 where the last inequality follows by $R-y<2t^{1/\alpha}$. If $R-y>2t^{1/\alpha}$, by \eqref{e1.56}  we have
  \begin{align*} 
 I&\leq C   \frac{t}{(R-y)^{ (2+\gamma) \alpha}}+C \frac{t^{1/\alpha-\gamma}}{(R-y)^{1+\alpha}}\nn\\
 &\leq   C\frac{t^\rho}{(R-y)^{ (1+\gamma+\rho) \alpha}}\leq C\frac{t^\rho}{(R-y+t^{1/\alpha})^{ (1+\gamma+\rho) \alpha}},
 \end{align*}
 where the second inequality uses $\rho \leq 1\wedge (1/\alpha-\gamma)$ and $R-y>2t^{1/\alpha}$. The last inequality follows from $R-y>2t^{1/\alpha}$. Combining the above two cases, we get \eqref{e1.56} implies that for any $-R< y< R$, 
 \begin{align}\label{be1.56}
 \int_0^t du \int p_{t-u}^R(y,x) \frac{1}{(R-x+u^{1/\alpha})^{(2+\gamma)\alpha}} dx \leq    C\frac{t^\rho}{(R-y+t^{1/\alpha})^{ (1+\gamma+\rho) \alpha}}.
 \end{align}
 By using the symmetry arguments as in \eqref{e0.78},  \eqref{e0.79}, and  \eqref{e0.01}, one can easily check that \eqref{be1.56} further gives that  for any $-R< y< R$, 
  \begin{align*} 
   &\int_0^t du \int p_{t-u}^R(y,x) \frac{1}{(R-|x|+u^{1/\alpha})^{(2+\gamma)\alpha}} dx  \leq C\frac{t^\rho}{(R-|y|+t^{1/\alpha})^{ (1+\gamma+\rho) \alpha}},
  \end{align*}
  thus giving Lemma \ref{l0.1}.\\
  
  It remains to prove \eqref{e1.56}. We will work with two cases: (1) $R-2t^{1/\alpha}<y<R$; (2) $-R<y<R-2t^{1/\alpha}$.\\

  \no   {\bf Case 1.}  For the case of $y\in (R-2t^{1/\alpha}, R)$, we have
   \begin{align*}
  t^{-1-\gamma}\leq  C \frac{t}{(R-y+t^{1/\alpha})^{ (2+\gamma) \alpha}}.
  \end{align*}
  Hence it suffices to show that
   \begin{align}\label{ae3.56}
I \leq   C t^{-1-\gamma}, \quad \forall y\in (R-2t^{1/\alpha}, R).
  \end{align}
  Notice that for $y\in (R-2t^{1/\alpha}, R)$,
   \begin{align}\label{e3.56}
 &\int_{t/2}^t du \int p_{t-u}^R(y,x) \frac{1}{(R-x+u^{1/\alpha})^{(2+\gamma)\alpha}} dx\nn\\
&\leq  \int_{t/2}^t  du \int  p_{t-u}^R(y,x)  u^{-(2+\gamma)}   dx \leq \int_{t/2}^t  u^{-(2+\gamma)} du \leq Ct^{-1-\gamma}.
 \end{align}
For each $y\in (R-2t^{1/\alpha}, R)$, define
 \begin{align} 
  J:=  \int_0^{t/2} du \int p_{t-u}^R(y,x)\frac{1}{(R-x+u^{1/\alpha})^{(2+\gamma)\alpha}} dx.
 \end{align}
 In view of \eqref{e3.56},  \eqref{ae3.56} follows if we show that
  \begin{align}\label{be8.4}
  J \leq   C t^{-1-\gamma},\quad     \forall y\in (R-2t^{1/\alpha}, R).
 \end{align}
 
Fix $y\in (R-2t^{1/\alpha}, R)$. Separate the integral in $J$ into three integrals over three regions: (i) $R-u^{1/\alpha}<x<R$; (ii) $y-(t-u)^{1/\alpha}<x<R-u^{1/\alpha}$; (iii)   $-R<x\leq y-(t-u)^{1/\alpha}$.  Denote each corresponding integral by $J^{(i)}, J^{(ii)}, J^{(iii)}$ respectively so that $J=J^{(i)}+ J^{(ii)}+ J^{(iii)}$.\\
  
      {\bf Part 1(i).} For the case of $R-u^{1/\alpha}<x<R$, we apply \eqref{5e6} to get
 \begin{align} \label{be8.1}
p_{t-u}^R(y,x) \leq C  \frac{(R-x)^{\alpha/2}}{\sqrt{t-u}} (t-u)^{-1/\alpha}\leq Ct^{-1/2-1/\alpha}(R-x)^{\alpha/2},
  \end{align}
  where the last inequality follows by $u\leq t/2$. Use the above to see that
   \begin{align*}
 J^{(i)}&\leq  Ct^{-1/2-1/\alpha} \int_0^{t/2}   du \int_{R-u^{1/\alpha}}^R (R-x)^{\alpha/2} \frac{1}{(R-x+u^{1/\alpha})^{(2+\gamma)\alpha}} dx \\
&\leq Ct^{-1/2-1/\alpha}  \int_0^{t/2} u^{-2-\gamma} du \int_{R-u^{1/\alpha}}^R (R-x)^{\alpha/2}  dx\\
&\leq Ct^{-1/2-1/\alpha}  \int_0^{t/2} u^{1/\alpha-3/2-\gamma} du.
 \end{align*}
 Notice that $\gamma<1/\alpha-1/2$ implies that $1/\alpha-3/2-\gamma>-1$. Thus from the above, we get the bound
  \begin{align*}
 J^{(i)}&\leq  Ct^{-1/2-1/\alpha}   (t/2)^{1/\alpha-1/2-\gamma}\leq Ct^{-1-\gamma}. 
 \end{align*}

{\bf Part 1(ii).}  Turning to the case of $y-(t-u)^{1/\alpha}<x<R-u^{1/\alpha}$, we use \eqref{be8.1} to see that
\begin{align*}
 J^{(ii)}&\leq  Ct^{-1/2-1/\alpha} \int_0^{t/2}   du \int_{y-(t-u)^{1/\alpha}}^{R-u^{1/\alpha}} (R-x)^{\alpha/2} \frac{1}{(R-x+u^{1/\alpha})^{(2+\gamma)\alpha}} dx \\
&\leq  Ct^{-1/2-1/\alpha} \int_0^{t/2}   du \int_{y-(t-u)^{1/\alpha}}^{R-u^{1/\alpha}}  \frac{1}{(R-x)^{(3/2+\gamma)\alpha}} dx \\
&=  Ct^{-1/2-1/\alpha} \int_0^{t/2}   du \int_{u^{1/\alpha}}^{R-y+(t-u)^{1/\alpha}} r^{-(3/2+\gamma)\alpha}dr.
 \end{align*}
 Use $R-y<2t^{1/\alpha}$ and $t-u\leq t$ to bound the above by
 \begin{align*}
 J^{(ii)}&\leq Ct^{-1/2-1/\alpha} \int_0^{t/2}   du \int_{u^{1/\alpha}}^{3t^{1/\alpha}} r^{-(3/2+\gamma)\alpha}dr.
 \end{align*}
Now we discuss the following three cases:\\
\no(1) If $(3/2+\gamma)\alpha>1$, then 
   \begin{align*}
J^{(ii)}&\leq  Ct^{-1/2-1/\alpha} \int_0^{t/2}     (u^{1/\alpha})^{1-(3/2+\gamma)\alpha} du\\
 &=  Ct^{-1/2-1/\alpha} \int_0^{t/2} u^{1/\alpha-3/2-\gamma} du\leq Ct^{-1-\gamma},
 \end{align*}
 where the last inequality follows by $\gamma<1/\alpha-1/2$.\\
\no (2)  If $(3/2+\gamma)\alpha<1$, then 
  \begin{align*}
J^{(ii)}&\leq  Ct^{-1/2-1/\alpha} \int_0^{t/2}     (3t^{1/\alpha})^{1-(3/2+\gamma)\alpha} du\\
 &= Ct^{-1/2-1/\alpha}  (t/2)  \cdot (3t^{1/\alpha})^{1-(3/2+\gamma)\alpha}\leq Ct^{-1-\gamma}.
 \end{align*}
\no (3) If $(3/2+\gamma)\alpha=1$, by \eqref{5e6} we use another bound on $p_{t-u}^R$: $$p_{t-u}^R(y,x) \leq  C(t-u)^{-1/\alpha}\leq Ct^{-1/\alpha},$$ where the last follows by $u\leq t/2$. Using this bound, we get
 \begin{align} \label{be89.1}
 J^{(ii)}&\leq  Ct^{-1/\alpha} \int_0^{t/2}   du \int_{y-(t-u)^{1/\alpha}}^{R-u^{1/\alpha}}  \frac{1}{(R-x)^{(2+\gamma)\alpha}}dx \nn\\
&=  Ct^{-1/\alpha} \int_0^{t/2}   du \int_{u^{1/\alpha}}^{R-y+(t-u)^{1/\alpha}} r^{-(2+\gamma)\alpha}dr.
 \end{align}
 Since $(2+\gamma)\alpha>(3/2+\gamma)\alpha=1$ by assumption, \eqref{be89.1} implies the following bound:
   \begin{align*}
 J^{(ii)}&\leq Ct^{-1/\alpha} \int_0^{t/2}    (u^{1/\alpha})^{1-(2+\gamma)\alpha} du.
 \end{align*}
 Noticing that $1/\alpha-2-\gamma=(3/2+\gamma)-2-\gamma=-1/2>-1$, we get
   \begin{align*}
 J^{(ii)}&\leq Ct^{-1/\alpha}  (t/2)^{1/\alpha-1-\gamma}\leq Ct^{-1-\gamma}.
 \end{align*}
 
 We conclude from the above three cases that $ J^{(ii)} \leq  Ct^{-1-\gamma}$ as required.\\

   {\bf Part 1(iii).} Considering the case of $-R\leq x\leq y-(t-u)^{1/\alpha}$, we get  
 \begin{align} \label{be8.2}
\frac{1}{(R-x+u^{1/\alpha})^{(2+\gamma)\alpha}}\leq   \frac{1}{ (R-y+(t-u)^{1/\alpha})^{(2+\gamma)\alpha}}\leq C  t^{-2-\gamma}.
  \end{align}
Recall \eqref{5e6} to see that
 \begin{align} \label{be8.3}
p_{t-u}^R(y,x) \leq C   \frac{t-u}{|x-y|^{1+\alpha}}.
\end{align}
Apply \eqref{be8.2} and \eqref{be8.3} to get
 \begin{align*}
 J^{(iii)}&\leq  Ct^{-2-\gamma} \int_0^{t/2}   (t-u) du  \int_{-R}^{y-(t-u)^{1/\alpha}}  \frac{1}{|x-y|^{1+\alpha}}  dx\\
 & =  Ct^{-2-\gamma} \int_0^{t/2}   (t-u) du  \int_{(t-u)^{1/\alpha}}^{y+R}  \frac{1}{r^{1+\alpha}}  dr\\
 &\leq  Ct^{-2-\gamma} \int_0^{t/2}   1 du \leq Ct^{-1-\gamma}.
 \end{align*}
The proof of \eqref{be8.4} is now complete.\\

  \no   {\bf Case 2.}  We turn to the proof of \eqref{e1.56} for $-R<y<R-2t^{1/\alpha}$, in which case
  \begin{align*}
&t^{1/\alpha-\gamma}(R-y)^{-1-\alpha}+t (R-y)^{-(2+\gamma)\alpha}\\
& \leq Ct^{1/\alpha-\gamma}(R-y+t^{1/\alpha})^{-1-\alpha}+Ct (R-y+t^{1/\alpha})^{-(2+\gamma)\alpha}.
  \end{align*}
Hence it suffices to show that
 \begin{align}\label{be8.10}
 I\leq C  \frac{t}{(R-y)^{ (2+\gamma) \alpha}}+C \frac{t^{1/\alpha-\gamma}}{(R-y)^{1+\alpha}}, \quad \forall -R<y<R-2t^{1/\alpha}.
 \end{align}
 Fix $-R<y<R-2t^{1/\alpha}$. Separate the integral in $I$ into five integrals over five regions: (i) $R-u^{1/\alpha}<x<R$; (ii) $(R+y)/2<x<R-u^{1/\alpha}$; (iii) $y+(t-u)^{1/\alpha}\leq x\leq (R+y)/2$; (iv) $y-(t-u)^{1/\alpha}\leq x\leq y+(t-u)^{1/\alpha}$; (v) $-R<x\leq y-(t-u)^{1/\alpha}$. Denote each corresponding integral by $I^{(i)}, I^{(ii)}, I^{(iii)}, I^{(iv)}, I^{(v)}$ respectively so that $I=I^{(i)}+I^{(ii)}+ I^{(iii)}+I^{(iv)}+I^{(v)}$.\\
  
      {\bf Part 2(i).} 
For the case of $x\in (R-u^{1/\alpha}, R)$, we get $x>(R+y)/2$ by $y<R-2t^{1/\alpha}$, thus giving $|y-x|\geq  (R-y)/2$. Use this and \eqref{5e6} to get
 \begin{align} \label{e1.79}
p_{t-u}^R(y,x) &\leq C  \frac{(R-x)^{\alpha/2}}{\sqrt{t-u}} \frac{t-u}{|x-y|^{1+\alpha}}\\
&\leq Ct^{1/2}(R-x)^{\alpha/2} \frac{1}{(R-y)^{1+\alpha}}, \quad  x\in (R-u^{1/\alpha}, R), y\in (-R, R-2t^{1/\alpha}).\nn
  \end{align}
By the above and $(R-x+u^{1/\alpha})^{-(2+\gamma)\alpha} \leq u^{-2-\gamma}$, we have
 \begin{align*}
 I^{(i)}&\leq Ct^{1/2}\frac{1}{(R-y)^{1+\alpha}} \int_0^t u^{-2-\gamma} du \int_{R-u^{1/\alpha}}^R (R-x)^{\alpha/2}  dx \\
&\leq Ct^{1/2} \frac{1}{(R-y)^{1+\alpha}}  \int_0^t u^{-2-\gamma} (u^{1/\alpha})^{1+\alpha/2} du \leq C\frac{1}{(R-y)^{1+\alpha}}  t^{1/\alpha-\gamma},
 \end{align*}
  where the last inequality follows by $\gamma<1/\alpha-1/2$.\\

{\bf Part 2(ii).} Turning to the case of $(R+y)/2<x<R-u^{1/\alpha}$, we still have $|y-x|\geq  (R-y)/2$. So we may use \eqref{e1.79} and that $(R-x+u^{1/\alpha})^{-(2+\gamma)\alpha} \leq (R-x)^{-(2+\gamma)\alpha}$ to see that in this case,
  \begin{align*}
 I^{(ii)}&\leq Ct^{1/2}\frac{1}{(R-y)^{1+\alpha}} \int_0^t   du \int_{(R+y)/2}^{R-u^{1/\alpha}}  (R-x)^{\alpha/2} \frac{1}{(R-x)^{(2+\gamma)\alpha}}  dx \\
&= Ct^{1/2} \frac{1}{(R-y)^{1+\alpha}}  \int_0^t  du \int_{u^{1/\alpha}}^{(R-y)/2} \frac{1}{r^{(3/2+\gamma)\alpha}} dr.
 \end{align*}
 Now we discuss the following three cases:\\
\no(1) If $(3/2+\gamma)\alpha>1$, then 
   \begin{align*}
 I^{(ii)}&\leq  Ct^{1/2} \frac{1}{(R-y)^{1+\alpha}}  \int_0^t     (u^{1/\alpha})^{1-(3/2+\gamma)\alpha} du \leq Ct^{1/\alpha-\gamma} \frac{1}{(R-y)^{1+\alpha}},
 \end{align*}
   where the last inequality follows by $\gamma<1/\alpha-1/2$.\\
 \no(2) If $(3/2+\gamma)\alpha<1$, then 
   \begin{align*}
 I^{(ii)}&\leq  Ct^{1/2} \frac{1}{(R-y)^{1+\alpha}}  \int_0^t    \Big(\frac{R-y}{2}\Big)^{1-(3/2+\gamma)\alpha} dr\\
 &\leq  Ct^{3/2} \frac{1}{(R-y)^{(5/2+\gamma)\alpha}}\leq C  \frac{t}{(R-y)^{(2+\gamma)\alpha}},
 \end{align*}
 where the last inequality follows by $R-y>2t^{1/\alpha}$.
 
 \no(3) If $(3/2+\gamma)\alpha=1$,  by \eqref{5e6} we use another bound on $p_{t-u}^R$: 
  \begin{align*}
 p_{t-u}^R(y,x) \leq  C \frac{t-u}{|y-x|^{1+\alpha}}\leq Ct {(R-y)^{-1-\alpha}},
  \end{align*}
where the last follows by $|y-x|\geq  (R-y)/2$. Using this bound, we get
 \begin{align} \label{e15.79}
 I^{(ii)}&\leq Ct \frac{1}{(R-y)^{1+\alpha}} \int_0^t   du \int_{(R+y)/2}^{R-u^{1/\alpha}}   \frac{1}{(R-x)^{(2+\gamma)\alpha}}  dx \nn\\
&= Ct  \frac{1}{(R-y)^{1+\alpha}}  \int_0^t   du \int_{u^{1/\alpha}}^{(R-y)/2} \frac{1}{r^{(2+\gamma)\alpha}} dr.
 \end{align}
 Since $(2+\gamma)\alpha>(3/2+\gamma)\alpha=1$, \eqref{e15.79} implies the following bound:
   \begin{align*}
 I^{(ii)}&\leq Ct  \frac{1}{(R-y)^{1+\alpha}}  \int_0^t   (u^{1/\alpha})^{1-(2+\gamma)\alpha} du.
  \end{align*}
 Noticing that $1/\alpha-2-\gamma=(3/2+\gamma)-2-\gamma= -1/2>-1$, we get
   \begin{align*}
 I^{(ii)}&\leq Ct  \frac{1}{(R-y)^{1+\alpha}}  t^{1/\alpha-1-\gamma}=Ct^{1/\alpha-\gamma} \frac{1}{(R-y)^{1+\alpha}}.
 \end{align*}
 
 We conclude from the above three cases that 
 \[
  I^{(ii)}\leq Ct  \frac{1}{(R-y)^{(2+\gamma)\alpha}}+ Ct^{1/\alpha-\gamma} \frac{1}{(R-y)^{1+\alpha}},
 \]
as required.\\
  
  {\bf Part 2(iii).}  For the case of $y+(t-u)^{1/\alpha}\leq x\leq (R+y)/2$, we get $R-x\geq (R-y)/2$ and so
 \begin{align} \label{e7.79}
\frac{1}{(R-x+u^{1/\alpha})^{(2+\gamma)\alpha}}\leq C   (R-y)^{-(2+\gamma)\alpha}.
  \end{align}
  By \eqref{5e6}, we get $p_{t-u}^R(y,x) \leq C   \frac{t-u}{|x-y|^{1+\alpha}}$. Use this and the above to see that 
 \begin{align*}
 I^{(iii)}&\leq  C (R-y)^{-(2+\gamma)\alpha} \int_0^t   (t-u) du  \int_{y+(t-u)^{1/\alpha}}^{(R+y)/2}  \frac{1}{|x-y|^{1+\alpha}}  dx\\
 & =  C (R-y)^{-(2+\gamma)\alpha} \int_0^t   (t-u) du  \int_{(t-u)^{1/\alpha}}^{(R-y)/2}  \frac{1}{r^{1+\alpha}}  dr\\
 &\leq  C (R-y)^{-(2+\gamma)\alpha} \int_0^t   1 \ du \leq Ct (R-y)^{-(2+\gamma)\alpha}.
 \end{align*}
 
 {\bf Part 2(iv).} Considering the case of ${y-(t-u)^{1/\alpha}}<x<{y+(t-u)^{1/\alpha}}$, we still have $R-x\geq (R-y)/2$ by $x< y+(t-u)^{1/\alpha}\leq (R+y)/2$, so use \eqref{e7.79} and $p_{t-u}^R(y,x) \leq C   (t-u)^{-1/\alpha}$ from \eqref{5e6} to get
 \begin{align*}
 I^{(iv)}&\leq  C (R-y)^{-(2+\gamma)\alpha} \int_0^t    du  \int_{y-(t-u)^{1/\alpha}}^{y+(t-u)^{1/\alpha}}  (t-u)^{-1/\alpha}   dx\\
  &=  C (R-y)^{-(2+\gamma)\alpha} \int_0^t  2\  du \leq Ct (R-y)^{-(2+\gamma)\alpha}.
 \end{align*}

   {\bf Part 2(v).} Finally for the case of ${-R}<x<{y-(t-u)^{1/\alpha}}$, we use \eqref{e7.79}  and $p_{t-u}^R(y,x) \leq C   \frac{t-u}{|x-y|^{1+\alpha}}$ from \eqref{5e6} to obtain 
 \begin{align*}
 I^{(v)}&\leq  C (R-y)^{-(2+\gamma)\alpha} \int_0^t   (t-u) du  \int_{-R}^{y-(t-u)^{1/\alpha}}  \frac{1}{|x-y|^{1+\alpha}}  dx\\
 & =  C (R-y)^{-(2+\gamma)\alpha} \int_0^t   (t-u) du  \int_{(t-u)^{1/\alpha}}^{y+R}  \frac{1}{r^{1+\alpha}}  dr\\
 &\leq  C (R-y)^{-(2+\gamma)\alpha} \int_0^t  1\  du \leq Ct (R-y)^{-(2+\gamma)\alpha}.
 \end{align*}
 The proof of \eqref{be8.10} is now complete by combining the above five parts.
 \end{proof}

 \end{document}